\documentclass{amsart}

\usepackage{amsmath}
\usepackage{amssymb}
\usepackage{amsfonts}
\usepackage{mathrsfs}
\usepackage{graphicx}
\usepackage{epsfig}
\usepackage{amsthm}
\theoremstyle{plain}

\newtheorem{proposition}{Proposition}[section]
\newtheorem{theorem}[proposition]{Theorem}
\newtheorem{lemma}[proposition]{Lemma}
{\theoremstyle{remark} \newtheorem{remark}[proposition]{Remark}}
\newtheorem{definition}[proposition]{Definition}

{\theoremstyle{remark} \newtheorem*{question}{Question} }
\newtheorem*{main theorem}{Main Theorem}

\def\proof{\par\noindent \textbf{Proof.} }
\newcommand{\proofof}[1]{\par\noindent \textbf{Proof of #1.}}
\def\finproof{$\hfill \square$\\}
\newcommand{\partn}[1]{{\smallskip \noindent \textbf{#1.}}}

\def\C{\mathbb{C}}
\def\D{\mathbb{D}}

\def\H{\mathbb{H}}

\def\R{\mathbb{R}}
\def\Z{\mathbb{Z}}

\def\T{\R / \Z}

\def\bfT{\mathbf{T}}

\def\bfx{\mathbf{x}}

\def\cB{\mathcal{B}}
\def\cC{\mathcal{C}}
\def\cD{\mathcal{D}}
\def\cF{\mathcal{F}}

\def\cO{\mathcal{O}}

\def\cR{\mathcal{R}}

\def\cK{\mathcal{K}}

\def\sA{\mathscr{A}}
\def\sF{\mathscr{F}}
\def\sG{\mathscr{G}}
\def\sR{\mathscr{R}}
\def\sT{\mathscr{T}}
\def\sU{\mathscr{U}}
\def\sV{\mathscr{V}}
\def\sW{\mathscr{W}}

\def\tB{{\widetilde{B}}}
\def\tC{{\widetilde{C}}}
\def\tD{{\widetilde{D}}}

\def\tf{{\widetilde{f}}}
\def\tg{{\widetilde{g}}}
\def\tF{{\widetilde{F}}}

\def\tq{{\widetilde{q}}}

\def\tX{{\widetilde{X}}}

\def\tgamma{{\widetilde{\gamma}}}
\def\tiota{{\widetilde{\iota}}}
\def\trho{{\widetilde{\rho}}}
\def\ttau{{\widetilde{\tau}}}
\def\tlambda{\widetilde{\lambda}}
\def\tPi{\widetilde{\Pi}}
\def\tSigma{\widetilde{\Sigma}}

\def\hF{{\widehat{F}}}
\def\hG{{\widehat{G}}}

\def\hp{{\widehat{p}}}
\def\hT{{\widehat{T}}}
\def\hU{{\widehat{U}}}
\def\hV{{\widehat{V}}}
\def\hW{{\widehat{W}}}

\def\hgamma{{\widehat{\gamma}}}
\def\hGamma{{\widehat{\Gamma}}}

\def\hOmega{{\widehat{\Omega}}}

\def\hcB{{\widehat{\cB}}}
\def\hcC{{\widehat{\cC}}}
\def\hcK{{\widehat{\cK}}}
\def\tcK{{\widetilde{\cK}}}
\def\tsU{{\widetilde{\sU}}}
\def\tsV{{\widetilde{\sV}}}

\def\ov{\overline}

\def\la{\lambda}

\DeclareMathOperator\dist{dist}

\DeclareMathOperator\interior{interior}
\DeclareMathOperator\Jac{Jac}
\DeclareMathOperator\length{length}

\def\sFla{\sF_\la}
\def\sGla{\sG_\la}
\def\sTla{\sT_\la}
\def\DH{H_\la}
\def\rad{\mathbf{r}}
\def\AF{A}
\def\AFla{\AF_\la}
\def\AFF{\AF_F}

\def\gcan{g_{\dag}}
\def\Gcan{G_{\dag}}
\def\hGcan{{\hG_{\dag}}}

\def\P1{\Pi_\la}
\def\tP1{\widetilde{\Pi}_\la}

\def\folss{{\cF}^{ss}}
\def\sing{{o}}
\def\tsing{{\widetilde{o}}}
\def\trapp{\mathbf{U}}
\def\pretrapp{\mathbf{V}}
\def\tpretrapp{{\widetilde{\pretrapp}}}
\def\linearpiece{\mathbf{W}}

\begin{document}
\title{Wild Lorenz like attractors}
\author{Rodrigo~Bam{\'o}n $\dag$,
Jan~Kiwi $\ddag$.
Juan~Rivera Letelier$\star$}
\thanks{$\dag$ Partially supported by Fondecyt N~1020711.}
\thanks{$\ddag$ Partially supported by Fondecyt N~1020711 and Fondecyt N~1010865.}
\thanks{$\star$ Partially supported by MeceSup PUC-UCH-0103 and MeceSup UCN-0202.}
\thanks{All authors were partially supported by IMPA and PROSUL}

\address{$\dag$Departamento de Matem\'aticas, Facultad de Ciencias,
Universidad de Chile, Casilla~653, Santiago, Chile.}  
\email{{\tt rbamon@uchile.cl}} 
\address{$\ddag$Facultad de Matem\'aticas,
Pontificia Universidad Cat\'olica de Chile, Casilla~306, Santiago~22,
Chile.} 
 \email{{\tt jkiwi@puc.cl}} 
\address{$\star$Departamento de
Matem\'aticas, Universidad Cat\'olica del Norte, Casilla~1280,
Antofagasta, Chile.} 
\email{{\tt juanrive@ucn.cl}}
\date{\today}
\maketitle

\begin{abstract}
  We give the first examples of flows which exhibit
  robust singular attractors containing a wild hyperbolic set (in the sense of Newhouse).
  A hyperbolic set is said to be wild, if it
  has tangencies between its stable and unstable manifolds, in a
  robust way.
  The only restriction on the ambient manifold is that its
  dimension should be at least $5$.
\end{abstract}


 \setcounter{tocdepth}{1} \tableofcontents


\section{Introduction}
In~\cite{L}, motivated by the equations for fluid convection in a two-dimensional layer heated
from below,  E.~Lorenz discovered a flow in $\R^3$
defined by an explicit and simple system of differential equations
that  exhibited extremely rich and robust dynamical behavior.
The combination of the works  \cite{ABS,GW,Tu} show  that
the Lorenz flow has a {\it singular attractor}: a transitive compact invariant
set that attracts a whole neighborhood of initial conditions and
that contains both  a singularity (stationary point) and regular orbits. Also, this
attractor is {\it robust}: every flow that is sufficiently close to the
initial one has a singular attractor nearby.
The topological, geometrical and ergodic properties of the Lorenz flow
have been intensively studied during the last two decades
(e.g., see~\cite{BS,BDV,LMP,M,S,V} and the references therein).
Moreover, the Lorenz flow plays a key role in a global picture for flows in three dimensions.
In fact, in~\cite{MPP} it is shown that every robust transitive set is either hyperbolic or
a singular attractor (or repellor) that shares many properties with the Lorenz attractor.

It is natural to ask whether higher dimensional manifolds support robust singular attractors with richer
dynamical behavior. This question has been addressed in~\cite{BPV} and~\cite{ShT}. While the singularity
of the Lorenz flow has $1$ expanding eigenvalue, Bonnatti, Pumari\~{n}o and Viana~\cite{BPV}
construct  robust singular attractors
in an $n \geq 4$ dimensional manifold where the singularity has $k = \max\{n-3, 2\} \geq 2$ expanding eigenvalues.
They also show that these attractors support a unique physical (SBR) measure with the no-holes property.
Shil'nikov and Turaev~\cite{ShT} construct a robust singular  {\em chain} transitive {\em chain} attracting set
 in $\R^4$  with the special feature that it contains a {\em wild hyperbolic set}, in the sense
of Newhouse~\cite{N}. Nevertheless, it is not known if this example has the key properties of being isolated and
transitive. A hyperbolic set is said
to be {\it wild}, if it has tangencies between its stable and
unstable manifolds, in a robust way.
By~\cite{MPP}, robust singular attractors in three dimensional
manifolds  are free of wild hyperbolic sets.
Moreover, the examples in~\cite{BPV} can not have a wild hyperbolic set since their dynamics reduces
to a  uniformly expanding endomorphism.

The purpose of this article is to give the first examples of
robust singular attractors containing a {\it wild hyperbolic
set}.  The only restriction on the
ambient manifold is that its dimension should be at least~5. In
fact, every flow on a manifold of dimension at least~5 having an
attracting periodic orbit, can be modified in an arbitrarily small
neighborhood of this orbit to have such a singular attractor.

\begin{main theorem}
Every manifold of dimension at least~5 admits a $C^1$ non-empty open set $
\cO $  such that every $C^2$ vector field $X$ in
$\cO$ exhibits a singular attractor $ \Lambda_X $ with the
following properties.
\begin{enumerate}
\item[1.] $ \Lambda_X $ contains a hyperbolic singularity of Morse
index ~$2$ of $X$. \item[2.] $ \Lambda_X $ contains a wild
hyperbolic set. \item[3.] There is a residual subset of $\cO$ that is dense in the $C^\infty$ topology,
such that if $X$ belongs to this residual set, then the set of
periodic orbits of $X$ of Morse index~$1$ and the set of periodic
orbits of $X$ of Morse index~$2$ are both dense in $ \Lambda_X $.
\end{enumerate}
\end{main theorem}

\begin{remark}
{ In fact 
we prove the stronger result.
that every $C^1$ vector field $X \in \cO$ has an attracting set $\Lambda_X$  satisfying properties
$1$, $2$, and $3$, see Theorem~\ref{t-dynamics of the flow} in Section~\ref{statements of results}. Although it is possible to push the ideas of this paper
to establish that $\Lambda_X$ is topologically transitive (and hence an attractor)
for all $X \in \cO$, it requires a fairly long and technical proof
(see~\cite{BKR}). }
\end{remark}

It is not clear to us what are the ergodic properties of the flows
introduced here. It would be very interesting to show, for
example, that these flows admit a physical or SBR measure.

\subsection{Strategy of the proof}
For the proof of the Main Theorem we proceed as follows.
Fix an integer $n \ge 5$. Then, for a given $\la \in (0,1)$ sufficiently close to $1$ we construct a family of
vector fields $\{ X_{\la, \mu} \}_{\la, \mu}$, where $\mu > 0$
takes values on a certain interval to be precised later.  These
vector fields are defined on a subset of the closed solid torus
$\bfT^n$ of dimension $n$. For some values of $\mu$ and a suitable
$C^1$-neighborhood $\cO$ of $X_{\la, \mu}$, we show that  the
conclusions of the Main Theorem are satisfied. The Main Theorem is
obtained by embedding the solid torus $\bfT^n$ on a given manifold
of dimension~$n$.

The vector field $X_{\lambda, \mu}$ will have a unique singularity
$o$. This singularity will be hyperbolic with eigenvalues $-\eta$,
$-\mu$ and $\sigma$ with multiplicities $n-3$, $1$ and $2$
respectively, where $0<\mu < \sigma < \eta$. The construction  is
such  that the dynamics of the flow of $X_{\la, \mu}$ reduces to the dynamics of the
map
\begin{eqnarray*}
F_{\lambda,\mu} : \C \setminus \{0\}  &  \to  &  \C \\
z  &  \mapsto  & ( 1 - \lambda + \lambda |z|^{\mu/\sigma} ) \left(
z / |z| \right)^2 + 1,
\end{eqnarray*}
in the same way as the dynamics of the geometric Lorenz attractor (see~\cite{ABS,GW})
reduces to the dynamics of a map defined on a punctured interval.
More precisely, the vector field $X_{\la,\mu}$ will induce a first return (or Poincar\'e) map $\hF_{\la,\mu}$ to a certain transversal section of
dimension $n - 1$. This first return map $\hF_{\lambda,\mu}$ will
have an invariant foliation of dimension $n - 3$  that is
uniformly contracted by $\hF_{\la,\mu}$. The map $F_{\la,\mu}$
defined above represents the leaf space  transformation of
$\hF_{\lambda,\mu}$; see Subsection~\ref{ss-flows} below for more
details.

The parameter $\mu$ will be taken in an interval of the form
$[\mu_0, \sigma]$, where $\mu_0 \in (0, \sigma)$ is sufficiently close to
$\sigma$. That is, $F_{\la,\mu}$ will be close to the map
\begin{eqnarray*}
F_\la :  \C \setminus \{0\}  &  \to  &  \C \\
z & \mapsto & (1 - \lambda + \lambda |z| ) \left( z / |z| \right)^2 +1.
\end{eqnarray*}
Note that the endomorphism $z \mapsto (1 - \lambda + \lambda |z| )
\left( z / |z| \right)^2$ acts as angle doubling on the argument
and as an affine contraction of factor $\lambda$ in the radial
direction. Thus $F_\lambda$ is closely related to the extensively
studied quadratic family $Q_c(z) = z^2 + c$, where the $|z|
\mapsto |z|^2$ action of $Q_0(z) = z^2$ is replaced by an affine
contraction of factor $\lambda$.

The point $0 \in \C$ represents a leaf contained in the stable
manifold of the  singularity $o$ of the vector field
$X_{\lambda,\mu}$. The map $F_{\lambda,\mu}$, which is not defined
at $z= 0$, `opens' the punctured plane $\C^* = \C \setminus \{0\}$
and maps it onto $\{ z \in \C \mid |z - 1| > 1 - \lambda \}$, in
a~$2$ to~$1$ fashion.

The dynamical properties of the maps $F_{\lambda,\mu}$, the
corresponding properties of the vector fields $X_{\lambda,\mu}$
and their relation to the Main Theorem are all explained in
Section~\ref{statements of results}.

\subsection{Notes and references}
To show that the attractor sets introduced here are robustly
transitive, we show that the stable and unstable manifolds of a
certain saddle fixed point are dense.
In our context, the  methods used
in~\cite{BPV} break down, since our leaf space transformation is not uniformly expanding.
A key step towards
showing that the unstable manifold is dense uses an argument
similar to that of~\cite{BD1}. In fact, our examples have some
kind of ``solenoidal blender''. The density of the stable manifold
is obtained by ad hoc arguments.

In~\cite{ShT} the
existence of a wild hyperbolic set is obtained, among other
ingredients, from the existence of  Newhouse phenomena of
persistence of tangencies.
Neither in~\cite{BD2} nor in our paper, the  persistence of
tangencies is obtained from the existence of  Newhouse
phenomena. Although the existence of a wild hyperbolic set implies
the co-existence of infinitely many periodic orbits of distinct
index~\cite{R}, we prove part~3 of the Main Theorem directly, {\it
independently of part~2}.

As mentioned above it is not clear to us what are the ergodic properties of the flows and endomorphisms
introduced here.  Note that when $\la \to 1$ the maps $F_\lambda$ converge to the map
$G(z) = |z| (z / |z|)^2 + 1$ in the $C^\infty$ topology. The
map~$G$ extends continuously to $\C$ and it preserves the Lebesgue
measure on $\C$.

\begin{question}
{ Is $G$ is ergodic with respect to the Lebesgue measure?}
\end{question}

\subsection{Acknowledgments}
We are thankful to M.~Viana and L.~D\'{\i}az, for useful comments.
Juan Rivera is grateful to C.~Morales for useful comments and for
pointing out the reference~\cite{ShT} and to B.~San Mart{\'i}n and
C.~V\'asquez for useful conversations.
Part of this paper was written while the authors visited IMPA, we would like express our gratitude to this institution for its hospitality and excellent working conditions.
Juan Rivera is also grateful to Pontificia Universidad Cat\'olica de Chile and to Universidad de Chile, for their hospitality.

\section{Statement of results}\label{statements of results}
In this section we outline our results and simultaneously describe
the structure of the paper. Recall that for the proof of the Main
Theorem we construct a family of
vector fields $X_{\la,\mu}$ so that the leaf space transformation of a
first return map to a cross section is close to
$$
\begin{array}{rccl}
  F_\la :&  \C \setminus \{0\} & \to & \C \\
            &  z                           & \mapsto & ( 1 - \lambda + \lambda |z| ) \left( {z}/{|z|} \right)^2  + 1.
\end{array}
$$
In Subsection~\ref{ss-the spaces} we introduce a topological space of $C^1$ maps $\sFla$ that contains the endomorphisms induced by $C^2$ vector fields which are close to $X_{\la,\mu}$ in the $C^1$ topology.
The main dynamical properties of maps in $\sFla$ close to $F_\la$ are summarized in Subsection~\ref{ss-attractor set}.
In Subsection~\ref{ss-flows} we describe the vector fields~$X_{\lambda,\mu}$ and in Subsection~\ref{ss-dynamics of the flow} we state a more precise version of the Main Theorem.
\subsection{The space $\sFla$}\label{ss-the spaces}

Throughout this paper the punctured complex plane $\C \setminus \{0\}$ will be denoted  by $\C^*$.
Note that the map $F_\lambda$ can be written as the composition $F_\lambda =
\gcan  \circ \tau_\la$ of the maps
$$ \begin{array}{rccl}
\tau_\la : & \C^* & \rightarrow & \T \times (1-\la, + \infty) \\
& z & \mapsto & (\tfrac{1}{2\pi} \arg(z), 1-\la + \la |z|)
\end{array}$$
and
$$
\begin{array}{rccl}
\gcan  : & \T \times (0, +\infty) & \rightarrow & \C \\
& (\theta, t) & \mapsto &  t \exp(4\pi i \theta)  + 1.
\end{array}
$$
Moreover note that  $\tau_\la$ is a diffeomorphism and that $\gcan$ is a local
diffeomorphism.

For $\la \in (0, 1)$ put
$$B_\lambda = \{ z \in \C \mid |z| \le 2(1 - \la)^{-1} \}$$
and
$$B_\lambda^* = B_\lambda \setminus \{ 0 \}.$$
It is easy to check that $F_\lambda (B_\lambda^*)$ is contained in
the interior of $B_\lambda$.

\begin{definition}\label{xla}
For a given $\lambda \in (0, 1)$ define the following spaces of $C^1$ maps.
\begin{enumerate}
\item[1.]
The space $\sTla$ of all homeomorphisms,
$$
\tau : B_\lambda^* \to \T \times (1 - \lambda, 2(1 - \lambda)^{-1} - 1 - \lambda],
$$
that extend to a diffeomorphism onto its image defined on a neighborhood of $B_\lambda^*$ in $\C^*$.
We endow $\sTla$ with the weak $C^1$ topology.
\item[2.]
The space $\sGla$ of all maps
$$
g : \T \times [1 - \lambda, 2(1 - \lambda)^{-1} -1 - \lambda] \to B_\lambda,
$$
having an extension to a neighborhood of $\T \times [1 - \lambda, 2(1 - \lambda)^{-1} - 1 - \lambda]$ in $\T \times \R$ that is a local diffeomorphism.
We endow $\sGla$ with the strong $C^1$ topology.
\item[3.]
The space $\sFla$ of all maps $F : B_\lambda^* \to B_\lambda$ that can be written as a composition $F = g \circ \tau$, with $\tau \in \sTla$ and $g \in \sGla$.
We endow $\sFla$ with the largest topology for which the composition map $(\tau, g) \mapsto g \circ \tau$, from $\sTla \times \sGla$ to $\sFla$, is continuous.
\end{enumerate}
\end{definition}
Clearly $F_\la \in \sFla$.
By definition, each map in $\sFla$ is a local diffeomorphism and hence it is an open map.

Notice that two maps in
$\sFla$ which  are close, need not be close with respect to
the uniform topology on~$B_\lambda^*$.
On the other hand, two maps in $\sF_\la$ that are close in the weak $C^1$ topology need not be close in $\sF_\la$.
In Lemma~\ref{l-composition map} below we show that the composition map from $\sTla \times \sGla$ to $\sFla$ is open.
It follows that for every $\tau \in \sTla$, $g \in \sG_\la$ and every $\tF \in \sFla$ that is close to $g \circ \tau$ in $\sFla$, there are $\ttau \in \sTla$ close to $\tau$ and $\tg \in \sGla$ close to $g$, so that $\tF = \tg \circ \ttau$.

The spaces $\sTla$ and $\sGla$ are Baire spaces as they both admit a complete metric.
Since the composition map from $\sTla \times \sGla$ to $\sFla$ is open (Lemma~\ref{l-composition map}), it follows that $\sFla$ is also a Baire space.

\begin{lemma}\label{l-composition map}
The composition map from $\sTla \times \sGla$ to $\sFla$ is open.
\end{lemma}
\proof
Denote the composition map by $\pi$.
We need to show that for every open set $\sW$ in $\sTla \times \sGla$ the set $\pi^{-1}(\pi(\sW))$ is open.
For that, let $(\ttau_0, \tg_0) \in \pi^{-1}(\pi(\sW))$ be given and let $(\tau_0, g_0) \in \sW$ be such that
\begin{equation}\label{e-0}
g_0 \circ \tau_0 = \pi((\tau_0, g_0)) = \pi((\ttau_0, \tg_0)) = \tg_0 \circ \ttau_0.
\end{equation}
Let $\sU$ (resp. $\sV$) be a neighborhood of $\tau_0$ (resp. $g_0$) in $\sTla$ (resp. $\sGla$) such that $\sU \times \sV \subset \sW$.

From~\eqref{e-0} it follows that the map $\ttau_0 \circ \tau_0^{-1}$ extends to a diffeormorphism onto its image $h_0$, defined on a neighborhood of $\T \times [1 - \lambda, 2(1 - \lambda)^{-1} - 1 - \lambda ]$.
By continuity we have $\tg_0 \circ h_0 = g_0$.
Put $\tsU = \{ h_0 \circ \tau \mid \tau \in \sU \}$ and $\tsV = \{ g \circ h_0^{-1} \mid g \in \sV \}$.
By construction the neighborhood $\tsU \times \tsV$ of $(\ttau_0, \tg_0)$ in $\sTla \times \sGla$ is such that $\pi(\tsU \times \tsV) = \pi(\sU \times \sV) \subset \pi(\sW)$.
This shows that $\pi^{-1}(\pi(\sW))$ is open in $\sTla \times \sGla$ and finishes the proof the lemma.
\finproof

\subsection{Dynamics of maps in $\sFla$ near $F_\lambda$}\label{ss-attractor set}
For a fixed $\lambda \in (0, 1)$ close to~$1$, we will be
interested on the dynamics of a given map $F$ in $\sFla$ near
$F_\lambda$.
Although a map $F$ in $\sFla$ is not defined at $z=0$, we will let
it act on the subsets of $B_\lambda$ by $F(U) = F(U \setminus \{ 0
\})$.
For a positive integer $m$ we will denote by $F^m$ the $m$-th iterate of
this action. Note that the image of an open subset of $B_\lambda$
under this action is again an open subset of $B_\lambda$.

Our main object of study will be the maximal invariant set of $F$ in~$B_\la$:
$$
\Omega_F = \cap_{m \ge 1} F^m(B_\lambda).
$$

It is easy to see that $p_\lambda = 1 + (1 - \la)^{-1}$ is a
saddle fixed point of $F_\lambda$. For a map $F \in \sFla$
close to $F_\lambda$, we denote by $p_F$ the saddle fixed point of
$F$ that is the continuation of $p_\lambda$. Of course, $p_F \in
\Omega_F$.

\begin{figure}
\begin{center}
\psfig{file=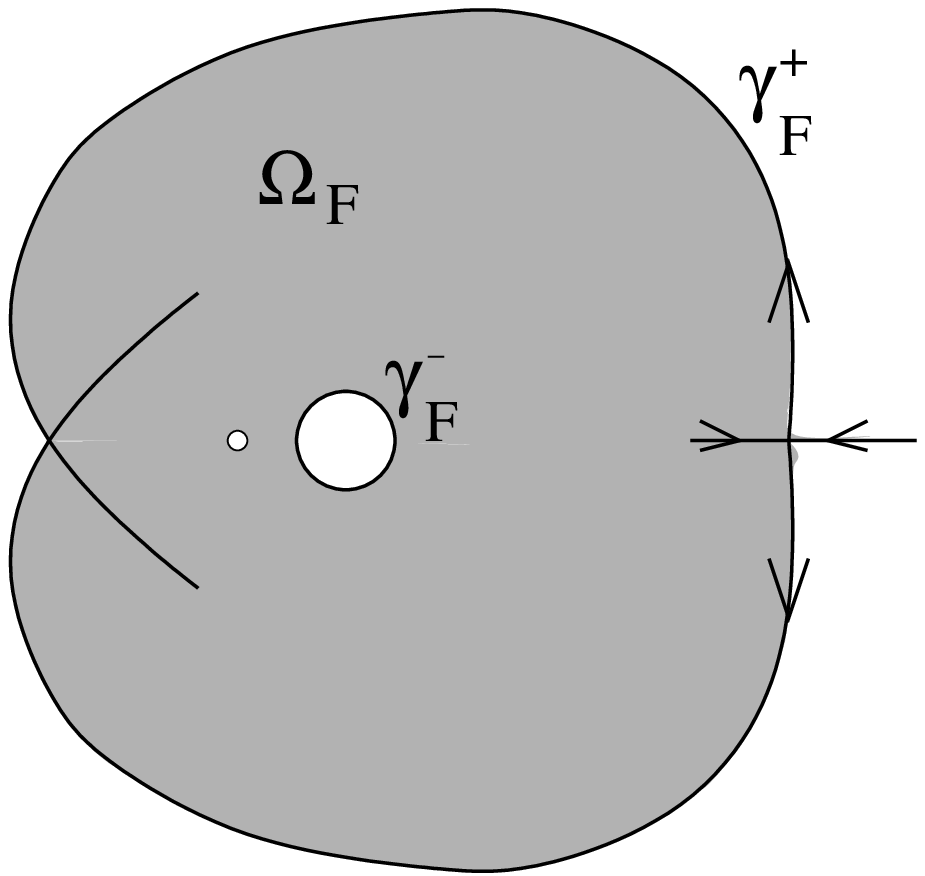,width=5cm}
\end{center}
\caption{The attractor $\Omega_F$} \label{omegaF-f}
\end{figure}

\begin{proposition}[The topology of $\Omega_F$]\label{p-maximal invariant}
For $\lambda \in (0, 1)$ sufficiently close to~$1$ and for $  F $
in  $\sFla$ sufficiently close to $F_\la$, the set
$\Omega_F$ is a topological  annulus that contains the origin  $z=
0 $ in its interior (see Figure~\ref{omegaF-f}).

The annulus $\Omega_F$ is neither open nor closed.  More
precisely, while the external  boundary $\gamma^+_F$ is contained
in $\Omega_F$, the  internal  boundary $\gamma^-_F$ is disjoint
from  $\Omega_F$.

Moreover, $\gamma^+_F$ is a Jordan curve contained in the unstable
manifold of $p_F$ and $\gamma^+_F$ is differentiable at every point with
the unique exception of a point where the unstable manifold of
$p_F$ has a transversal  self-intersection. The internal  boundary
$\gamma^-_F$ is a Jordan curve.
\end{proposition}

Section~\ref{technical results} contains the proof of this
proposition and some technical results which will be needed in the
rest of the paper.

In Section~\ref{fundamental annulus} we study some dynamical
features of maps close to $F_\la$ such as the existence of a
fundamental annulus with some important covering properties.

\medskip

When $\lambda \in (2^{- 1/2}, 1)$, one of the important features of the
map $F_\lambda$ is that its Jacobian,
$$
\Jac (F_\lambda)(z) = 2\lambda (\lambda + (1 - \lambda)|z|^{-1}),
$$
is everywhere larger than a constant larger than~$1$. We say that
a map $F$ in $\sFla$ is {\it area expanding}, if there is a
constant $\kappa > 1$ such that for every point $z \in
B_\lambda^*$ the Jacobian of $F$ at $z$ is at least $\kappa$.
A map in $\sFla$ near $F_\lambda$ need not be area
expanding, as there is no control on its behavior near
$z=0$.

In Section~\ref{dinamica atractor} we prove that for $\la \in (0,
1)$ sufficiently close to~$1$ the map $F_\la$ is robustly
transitive on the maximal invariant set $\Omega_{F_\la}$, in the
subspace of~$\sFla$ of area expanding maps. More precisely, we
prove the following result.

\begin{theorem}[Robust transitivity]\label{t-topological dynamics}
For $\lambda \in (0, 1)$ sufficiently close to~$1$ there exists a
neighborhood $\sU$ of $F_\la$ in $\sFla$ such that each area
expanding map $F$ in $\sU$ is topologically mixing on $\Omega_F$.
\end{theorem}

To prove the previous theorem we show that the stable
manifold of $p_F$ is dense in $B_\la$ and that the unstable
manifold of $p_F$ is dense in  $\Omega_F$. By a standard
argument we conclude that $F$ is topologically mixing on
$\Omega_F$.

In Section~\ref{s-periodic points} we show the following result about periodic points of maps in $\sFla$.
Given an open set $\sU$ of $\sGla$ and $\tau \in \sTla$ put
$$
\sU_\tau = \{ g \in \sU \mid F = g \circ \tau \text{ is area expanding}\}.
$$
Note that $\sU_\tau$ is an open subset of $\sGla$, that might be empty.
For example, if $\Jac(\tau)(z) \to 0$ as $z \to 0$ then the set $\sU_\tau$ is empty.

\begin{proposition}[Periodic points] \label{p-periodic points}
There is a neighborhood $\sU$ of $\gcan$ in $\sGla$ such that for every $\tau \in \sTla$ there is a residual subset $\sR_\tau$ of $\sU_\tau$ with the following property.
For every $g \in \sR_\tau$ the set of periodic sources and the set of periodic saddles of $F = g \circ \tau$ are both dense in $\Omega_F$.
\end{proposition}

In order to state the results concerning robust homoclinic
tangencies, we need to introduce some notation.
For $ \la \in (0, 1)$ close enough to $ 1 $ and $ F \in \sFla $ close enough to $
F_\la $, we will define two subsets of $\Omega_F$: a fundamental annulus $ \AFF $ (Section~\ref{fundamental annulus}) and a domain  $\DH$  containing $p_F$, where $ F $ is uniformly hyperbolic (Subsection~\ref{domain of hyperbolicity}).
The set $\AFF$ is such that
$$
\AFF \subset F(\AFF) \subset \DH,
$$
so the set
$$
\Gamma_F = \left\{ \{ z_j \}_{j \ge 0} \mid z_j \in \DH \text{ and
} F(z_{j + 1}) = z_j, \text{ for } j \ge 0 \right\},
$$
is non empty.  Moreover,  every infinite backward orbit $ \{ z_j
\}_{j \ge 0}$ in $\Gamma_F$ defines a local  unstable manifold
passing through $z_0$ (Subsection~\ref{ss-local unstable
manifolds}). In Section~\ref{robust tangencies} we prove the
following result.

\begin{theorem} [Homoclinic tangencies and wild hyperbolic sets]\label{t-wild for endos} For $\lambda \in (0, 1)$ close to~$1$ and every area expanding map $F$ in $\sFla$ close to $F_\la$ the following properties hold.
\begin{enumerate}
\item[1.]
The set
$$
W_F = \{z \in \DH \mid F^m(z) \in \DH \;
\; \text{ for every } \; m \ge 1\}
$$
is an uniformly hyperbolic and forward invariant set for $F$.
Moreover, the local unstable manifold of an infinite backward orbit in $\Gamma_F$ starting at a point in $\AFF$, is contained in the unstable manifold of some infinite backward orbit contained in~$W_F$.
\item[2.]
There is an arc $\tgamma_\la$ of the stable manifold of $p_\la$, such that if $\tgamma_F$ is an arc of the stable manifold of the fixed point $p_F$ of $F$ that is $C^1$ close to $\tgamma_\la$, then $\tgamma_F$ is tangent to the unstable manifold of an infinite backward orbit in $\Gamma_F$.
\end{enumerate}
As the saddle fixed point $p_F$ of $F$ is contained in $W_F$, it follows that the hyperbolic set $W_F$ of $F$ is wild.
\end{theorem}

Each of the  unstable manifolds associated to an infinite
backward orbit in $\Gamma_F$ is approximated in the $C^1$ topology
by a piece of the unstable manifold of $p_F$
(Subsection~\ref{ss-unstable density}). So by an arbitrarily small
perturbation  it is possible to create a homoclinic tangency of
the fixed point $p_F$.

\subsection{The vector fields $X_{\lambda,\mu}$}\label{ss-flows}
For an integer $k \ge 0$ denote by $D^k$ the closed unit ball of $\R^k$.
Fix an integer $n \ge 5$ and let $\bfT^n = (\R/\Z) \times D^{n-1}$
be the closed solid torus of dimension~$n$.

For the proof of the Main Theorem we fix constants $\eta > \sigma > 0$.
Given $\lambda \in (0,1)$ sufficiently close to $1$ we
construct a one  parameter family $\{ X_{\lambda,\mu} \mid \mu \in
(\mu_0,\sigma] \}$ of vector fields, where $\mu_0 \in (0,\sigma)$ is close to $\sigma$, defined on an open subset~$\trapp$ of $\bfT^n$.
For $\mu \in (\mu_0, \sigma]$ the vector field $X_{\la, \mu}$ will have a hyperbolic singularity~$o = o_{\la, \mu}$ with eigenvalues $-\mu$, $\sigma$ and $-\eta$ of multiplicities~$1$, $2$ and $n - 3$, respectively.

The vector fields $X_{\la,\mu}$ are closely related to
the map $F_\la$ and the space $\sFla$.
More precisely, for $\mu$ close to
$\sigma$  the map
$$
\begin{array}{rccl}
  F_{\la, \mu} : & B^*_\la & \rightarrow & B_\la \\
                       &  z          &  \mapsto    &  (1 - \la +\la |z|^{\mu/\sigma}) \left( z/|z| \right)^2  +1.
\end{array}
$$
is an area expanding endomorphism in $\sFla$ that is close to $F_\la$.
Next  we consider the skew product map
$$
\begin{array}{rccl}
  \hF_{\la,\mu} : & B^*_\la \times \overline{\D} \times D^{n-5} & \to         & B_\la \times \overline{\D} \times D^{n-5} \\
                         & (z,w,v)   &\mapsto &   \left( F_{\la,\mu}(z),  \dfrac{z}{ 2 |z|}
 + \beta |z|^{\eta/\sigma} \dfrac{|z|}{z} w,  \beta |z|^{\eta/\sigma} v \right),
\end{array}
$$
where $\beta>0$ is an arbitrarily small constant so that
$\hF_{\la,\mu}$ is well defined and injective.
The vector field $X_{\la, \mu}$ will be constructed in such a way that the first return map to a certain cross section $\Sigma^u$ of the flow of $X_{\la, \mu}$, parametrized by $B_\la \times \overline{\D} \times D^{n-5}$, is given by~$\hF_{\la, \mu}$.

Note that the base dynamics of $\hF_{\la, \mu}$ is given by $F_{\la,\mu}$.
More precisely, the fibers of the linear projection $\P1 : B_\la \times \overline{\D} \times D^{n-5}  \rightarrow B_\la$ onto the first coordinate form an
invariant (strong stable) foliation for the dynamics of $\hF_{\la,\mu}$ and the action induced by $\hF_{\la,\mu}$ on the leaf space of this foliation is exactly $F_{\la,\mu}$.
The condition $ \mu \in (0, \sigma)$ guarantees that for every vector field $X$ that is close to $X_{\la, \mu}$, the first return map $\hF_X$ to $\Sigma^u$ for the flow of $X$, will have an invariant strong stable foliation close to the one formed by the fibers of $\P1$.
When $X$ is of class $C^2$ the corresponding leaf space transformation lies in the space $\sFla$, is area expanding and is close to $F_{\la, \mu}$, and hence to $F_\la$ (Lemma~\ref{l-leaf space transformation}).
So we will be able to apply the results for singular endomorphisms, described in Subsection~\ref{ss-attractor set}, to study the flow of the vector fields $X_{\la,\mu}$ and its perturbations.

\smallskip
Although $\hF_{\la,\mu}$ depends on $\beta >0$, for simplicity of
the notation we omit it.
In theorems~\ref{embedding in R5} and~\ref{t-dynamics of the flow} below we choose $\beta$  to be
conveniently small.

\begin{theorem}\label{embedding in R5}
Fix an integer $n \ge 5$ and $\eta> \sigma >0$. Then for each $\la
\in (0,1)$ sufficiently close to $1$, there exists $\mu_0 \in
(0,\sigma)$ and a smooth one parameter family of smooth vector
fields  $\{ X_{\la,\mu} \mid \mu \in (\mu_0,\sigma] \}$ defined on
an open set $\trapp$ of $\bfT^n$, such that for all $\mu \in
(\mu_0,\sigma]$ the following hold:

\begin{enumerate}
\item[1.]
The boundary of the open set $\trapp  \subset \bfT^n$ is a manifold of dimension~$n - 1$ that is contained in the interior of $\bfT^n$.
For every $\mu \in [\mu_0, \sigma]$ the vector field $X_{\la, \mu}$ extends to a smooth vector field defined on a neighborhood of the closure of $\trapp$, in such a way that on the boundary of $\trapp$ this vector field points inward.
Moreover this extension has a unique singularity $o=o_{\la,\mu}$.
The singularity $o$ is contained in $\trapp$ and is hyperbolic with eigenvalues $-\mu$, $\sigma$ and $-\eta$ of multiplicities $1$, $2$ and $n-3$, respectively.

\item[2.]  There exists a  codimension~$1$  submanifold $\Sigma^u
\cong B_\la \times \overline{\D} \times D^{n-5} \subset \trapp$
transversal to the flow of $X_{\la,\mu}$  so that  every forward
orbit of the flow of  $X_{\la,\mu}$ in $\trapp$ intersects $\Sigma^u$
or is contained in a local stable manifold $W^s_{loc}(o)$ of $o$.
The intersection  of this local stable manifold  with $\Sigma^u$
is $\{ 0 \} \times \overline{\D} \times D^{n-5}$.

\item[3.] The  first return map to $\Sigma^u$ is given by
$$\hF_{\la,\mu} :  B^*_\la \times \overline{\D} \times D^{n-5}
\rightarrow B_\la \times \overline{\D} \times D^{n-5}. $$
\end{enumerate}
\end{theorem}
\subsection{Dynamics of vector fields near $X_{\la, \mu}$}\label{ss-dynamics of the flow}
Part~$1$ of Theorem~\ref{embedding in R5} implies that for every vector field $X$ that is sufficiently close to $X_{\la, \mu}$ and for every $t > 0$ we have $\ov{X^t(\trapp)} \subset \trapp$.
We denote by,
$$
\Lambda_X = \cap_{t  > 0} \, \overline{X^t (\trapp)},
$$
the maximal invariant set of the flow of~$X$ in~$\trapp$.
Moreover we denote by~$\sing_X$ the singularity of~$X$ that is the continuation of the hyperbolic singularity~$\sing$ of $X_{\la, \mu}$.
Clearly we have $\sing_X \in \Lambda_X$.
We will show that for $\lambda \in (0, 1)$ sufficiently close to~$1$ and $\mu \in
(\mu_0,\sigma)$ sufficiently close to $\sigma$, there is a $C^1$
neighborhood of $X_{\la,\mu}$ so that statements 1 through 3 of
the Main Theorem hold for all vector fields $X$ in~$\cO$.
More precisely we prove the following result, which implies the Main Theorem.
\begin{theorem}\label{t-dynamics of the flow}
For each $ \la  \in (0,1)$ sufficiently close to $1$ and each $\mu \in (0,\sigma)$ sufficiently close to $\sigma$ there exist a $C^1$-neighborhood $\cO$ of $X_{\la, \mu}$, such that for each $X \in \cO$ the maximal invariant set $\Lambda_X$ of the flow of $X$ in $\trapp$ satisfies the following properties.
\begin{enumerate}
\item[1.] If $X \in \cO$ is of class $C^2$, then the flow of $X$ restricted to $\Lambda_X$ is topologically mixing.
Moreover, $\sing_X \in \Lambda_X$ and~$\sing_X$ is the unique singularity of~$X$ contained in $\trapp$.
In particular $\Lambda_X$ is a robustly transitive attractor set in the $C^2$ topology.

\item[2.] The attractor set $ \Lambda_X $ contains an invariant
hyperbolic set $W_X $ having a tangency between its stable and
unstable foliations.

\item[3.] For vector fields in a residual subset of $ \cO, $ the
sets of periodic orbits of Morse index~1,
and Morse index~2 respectively,  are both dense in~$\Lambda_X$.
\end{enumerate}
\end{theorem}

The hypothesis in part~$1$, that $X \in \cO$ is of class $C^2$, is unnecessary: the statement still holds for vector fields of class~$C^1$.
The proof of this fact is more involved and is done in detail in~\cite{BKR}.

Observe that the singularity $o$ of $X_{\la, \mu}$ has a double
real unstable  eigenvalue~$\sigma$. Thus, the singularity $o_X$ of
vector fields near $X_{\la, \mu}$ may have real or complex
unstable eigenvalues.

\section{Domain of uniform hyperbolicity and the topology of the attractor}\label{technical results}

In this section we construct  stable and unstable
cone fields (subsections ~\ref{ss-stable cone field}
and~\ref{ss-unstable cone field}) which are invariant under $F_\la$ on  certain
domain $\DH$, that we call the {\it hyperbolicity domain}
(Subsection~\ref{domain of hyperbolicity}). After studying a
saddle fixed point $p_\la$ of $F_\la$ in the hyperbolicity domain
$\DH$, we prove Proposition~\ref{p-maximal invariant}
in Subsection~\ref{ss-fixed point}.

\

Notice that for every  $a, b \in \R $ we have

\begin{equation}\label{derivative}
\begin{split}
   D_{z_0}F_{\la}(z_0(a + ib))  & = \left(\frac{z_0}{|z_0|}\right)^2
    (a \la  |z_0| + 2 b i (1 - \la + \la|z_0|))  \\
     & = (F_{\la}(z_0) - 1) \left({\frac{a \la |z_0|}{1 - \la +
     \la|z_0|}} + 2 b i  \right).
\end{split}
\end{equation}

\subsection{Stable cone field}\label{ss-stable cone field}

Set $\varepsilon_0 = \tfrac{1}{2} \sqrt{1 - \la}$ and let $\{
\cC_z \}_{z \in \C^*}$ be the cone field defined by $$ \cC_z = \{ \rho z
(1 + i\varepsilon) \mid \rho \ge 0, |\varepsilon| \leq
\varepsilon_0\}.$$

\begin{lemma}\label{l-stable cone field}
\

\begin{enumerate}
\item[1.]
For every $z_0 \in \C$ such that $|z_0| \geq 1 $ and every $v \in
\cC_{z_0}$ we have $$ |D_{z_0}F_{\la}(v)|^2 \le \tlambda \cdot
|v|^2, \ {\rm where} \ \tlambda = {\la^2 + 4\varepsilon^2_0 \over
1 + \varepsilon^2_0} \in (0, 1).$$
\item[2.]
For $ z_0 \in \C $ satisfying $1 - \la + \la |z_0| \geq 4 /
\sqrt{1 - \la} $ we have $$\cC_{F_{\la}(z_0)} \subseteq
\interior(D_{z_0} F_{\la}(\cC_{z_0})) \cup \{ 0 \}.$$
\end{enumerate}
\end{lemma}
\proof
\

\partn{1}
Letting $v = z_0 (1 + i \varepsilon) $ with $|\varepsilon| \leq \varepsilon_0 $, we have

\begin{multline*}
|D_{z_0}F_{\la}(v)|^2 = \left| \la |z_0| + 2  \varepsilon i (1 - \la + \la|z_0|) \right|^2
= \\ =
(\la  |z_0|)^2 + (2  \varepsilon (1 - \la + \la|z_0|))^2.
\end{multline*}
As by hypothesis $|z_0| \ge 1$, we have
$$
|D_{z_0}F_{\la}(v)|^2 \le |z_0|^2 (\la^2 + 4 \varepsilon^2)
\le
|v|^2 {\la^2 + 4\varepsilon^2 \over 1 + \varepsilon^2} \le |v|^2 {\la^2 + 4\varepsilon^2_0 \over 1 + \varepsilon^2_0}.
$$

\partn{2}
Note that $$ |F_{\la}(z_0) - 1| = 1 - \la + \la|z_0| \geq
4 / \sqrt{1 - \la} = 2 \varepsilon_0^{-1}.$$

On the other hand, $$D_{z_0} F_{\la}(\cC_{z_0}) = \left\{
(F_{\la}(z_0) - 1)(a + ib) \in \C \mid |b / a | \leq 2
\varepsilon_0 (1 - \la + \la|z_0|) / (\la |z_0|) \right\}.$$

Let $ \widetilde{a}, \widetilde{b} \in \R $ be so that $
|\widetilde{b} / \widetilde{a}| \leq \varepsilon_0.$ Then
$F_{\la}(z_0)(\widetilde{a} + i \widetilde{b}) \in
\cC_{F_{\la}(z_0)}.$ Define $ a, b \in \R$ by $ a + ib =
F_{\la}(z_0) (F_{\la}(z_0) - 1)^{-1} ( \widetilde{a} + i
\widetilde{b}).$ We want to prove that
$$\left| \frac{b}{a} \right| < 2 \varepsilon_0 \frac{1 - \la + \la|z_0|}{\la |z_0|}.$$

Set $ F_{\la}(z_0)(F_{\la}(z_0) - 1)^{-1} = 1 + \eta_0 + i \eta_1$ and notice that we have $$ |\eta_0|, |\eta_1| \leq |F_{\la}(z_0) - 1|^{-1} \le 2 \varepsilon_0 < \tfrac{1}{4},$$ so  $|\eta_1 / (1 + \eta_0)| \leq \tfrac{2}{3} \varepsilon_0$.
On the other hand, $ a = \widetilde{a} (1 + \eta_0) - \eta_1 \widetilde{b} $, $b
= \widetilde{b} (1 + \eta_0) + \eta_1 \widetilde{a}$ and
$$ \frac{b}{a} = \frac{\widetilde{b} / \widetilde{a} (1 + \eta_0) + \eta_1}{1 + \eta_0 - \eta_1 \widetilde{b} / \widetilde{a}}
=
\frac{\widetilde{b} / \widetilde{a} + \eta_1 /(1 + \eta_0)}{1 - (\eta_1 / (1 + \eta_0)) (\widetilde{b} / \widetilde{a})}.$$

Therefore,

\begin{equation*}
\left| \frac{b}{a} \right| \le
\frac{\varepsilon_0 + \tfrac{2}{3} \varepsilon_0}{1 - \tfrac{2}{3} \varepsilon_0^2}  \le 2 \varepsilon_0
<
\frac{2 \varepsilon_0 (1 - \la + \la|z_0|)}{\la |z_0|}.
\end{equation*}

\finproof


\subsection{Unstable cone fields}\label{ss-unstable cone field}
Let $\{ \cK(z) \}_{z \in \C^*}, \; \{\cK^-(z) \}_{z \in \C^*}$ and
$\{ \tcK(z) \}_{z \in \C^*}$ be the cone fields defined by
\begin{eqnarray*}
\cK(z_0) & = & \{ z_0 \rho (i + \varepsilon) \mid \rho \ge 0, \
|\varepsilon| \le 1/3 \},
\\
\cK^-(z_0) & = & \{z_0 \rho (i + \varepsilon) \mid \rho \ge 0, \
\varepsilon \in [-1/3, 0] \},
\end{eqnarray*}
Note that for every $z \in \C^*$ we have $\cK(z) \subset \tcK(z)$.

\begin{lemma}[Unstable Cone fields]\label{l-unstable cone field}

\

\begin{enumerate}
\item[1.] For every $z_0 \in \C^*$ and every $v \in \tcK(z_0)$ we
have
$$
| D_{z_0}F_\lambda (v) | \ge \lambda (5/2)^{1/2} | v |.
$$
\item[2.] For every $z_0 \in \C^*$ such that $1 - \lambda +
\lambda |z_0| \ge 20 $ we have
$$
D_{z_0}F_\lambda (\cK(z_0)) \subset {\rm
interior}(\cK(F_\lambda(z_0))) \cup \{ 0 \} \text{ and}
$$
$$
\cK(F_\lambda(z_0)) \subset \interior (D_{z_0} F_\la(\tcK(z_0)))
\cup \{ 0 \}.
$$

\item[3.] For every $z_0 \in \C^*$ such that $1 - \lambda +
\lambda |z_0| \ge 20, \, \Re(z_0) \geq 0 $ and $ \Im(z_0) > 0$, we
have
$$
D_{z_0}F_\lambda (\cK^-(z_0)) \subset {\rm
interior}(\cK^-(F_\lambda(z_0))) \cup \{ 0 \}.
$$
\end{enumerate}
\end{lemma}

 \proof

\noindent {\bf 1.} Consider $v = z_0 (i + \varepsilon)$ with $|\varepsilon| \le 1$, so that $v \in \tcK(z_0)$.
Setting $a = \rho \varepsilon$ and $b = \rho$ in
(\ref{derivative}), we have
$$
|D_{z_0}F_\lambda(v)|^2 = (\rho|z_0|)^2 \left| \lambda \varepsilon
+ 2i((1 - \lambda)/|z_0| + \lambda) \right|^2 \ge |v|^2 \lambda^2
\frac{\varepsilon^2 + 4}{\varepsilon^2 + 1}.
$$

\noindent {\bf 2.} Note that $|F_\la(z_0) - 1| = 1 - \lambda +
\lambda |z_0| \ge 20$, so $|F_\la(z_0)| \ge 19$. Given $\rho, \varepsilon \in \R$ such that $\rho > 0, |\varepsilon| \le
1/3$, define $\rho' \ge 0$ and $\varepsilon' \in \R$ by
$$
D_{z_0} F_\lambda(\rho z_0 (i + \varepsilon)) = F_\la(z_0) \rho'(i
+ \varepsilon').
$$
It is enough prove that $|\varepsilon' - \varepsilon/2| < 1/6$.

Define $\rho'', \rho_0 \ge 0$ and $\varepsilon'', \varepsilon_0
\in \R$ by
$$
D_{z_0} F_\lambda(\rho z_0 (i + \varepsilon)) = (F_\la(z_0)  -
1)\rho''(i + \varepsilon''),
$$
and $F_\la(z_0) - 1 = F_\la(z_0)\rho_0(1 + i\varepsilon_0)$.

Then $\varepsilon' = \frac{\varepsilon'' - \varepsilon_0}{1 +
\varepsilon'' \varepsilon_0}$. From ($\ref{derivative}$) we have
$\varepsilon'' = (\varepsilon / 2) \lambda |z_0| (1 - \lambda +
\lambda |z_0|)^{-1}$. So $|\varepsilon''| < |\varepsilon /2| \le
1/2$ and
$$
|\varepsilon'' - \varepsilon/2| < |\varepsilon / 2| (1 - \lambda +
\lambda |z_0|)^{-1} \le 1/40.
$$
On the other hand $|\rho_0(1 + i\varepsilon_0) -1| =
|F_\la(z_0)|^{-1} \le 1/19$. Hence $|\rho_0 - 1|$, $|\rho_0
\varepsilon_0| \le 1/19$ and $|\varepsilon_0| \le 1/18$. Since
$|\varepsilon''| < 1/2$, we have
$$
|\varepsilon' - \varepsilon''| = \left| \varepsilon_0 \frac{1 +
(\varepsilon'')^2}{1 + \varepsilon'' \varepsilon_0} \right| \le
2|\varepsilon_0| \le 1/9.
$$
Therefore $|\varepsilon' - \varepsilon/2| \le 1/40 + 1/9 < 1/6$.

 \noindent {\bf 3.}
Let $z_0 \in \C^*$ be such that $1 - \la + \la|z_0| \ge 20$ and let $\rho' \ge 0$ and $\varepsilon' \in \R$ be such
that
$$
D_{z_0} F_\la(i z_0) = F_\la (z_0) \rho' (i + \varepsilon').
$$
By  part 2 is enough to prove that when $\Re z_0 \ge 0$ and $\Im
z_0 > 0$ we have $\varepsilon' < 0$. Note that when $\Re z_0 \ge
0$ and $\Im z_0 > 0$ we have $\Im F_\la(z_0) > 0$.

By ($\ref{derivative}$) we have that $D_{z_0}F_\la (i z_0) =
(F_\la(z_0) - 1) 2 i$. So
$$
\varepsilon'
=
- \Im \left( \frac{F_\la(z_0) - 1}{F_\la(z_0)} \right)
\left/
\Re \left( \frac{F_\la(z_0) - 1}{F_\la(z_0)} \right)
\right.
=
- \frac{\Im(F_\la(z_0))}{|F_\la(z_0)|^2 - \Re F_\la(z_0)} < 0.
$$
\finproof


\subsection{Domain of uniform hyperbolicity}\label{domain of
hyperbolicity}

For $ \la \in (0, 1)$ sufficiently close to ~$1$, the subset $\DH$
of $\C$ defined by
$$\DH = \left\{ z \in \C \mid |z| > \la^{-1} \left( \la - 1 + 4 / \sqrt{1 - \la} \right) \right\},$$ will be called the {\it domain of hyperbolicity}.
We will say that a differentiable arc $\gamma$ in $\C^*$ is {\it
quasi-radial} (resp. {\it quasi-angular}), if for every point $z$
of $\gamma$ the tangent space of $\gamma$ at $z$ is contained in
$\cC_z$ (resp. $\cK_z$).

For every map $F$ in $\sFla$ sufficiently close to $F_\la$, the
cone field $\cC$ (resp. $\cK$) is an stable (resp. unstable) and
invariant cone field of $F : \DH \to F(\DH)$. For such $F$, the image of a quasi-angular arc contained in $\DH$ is a quasi-angular arc whose length is larger than the original length by a definite factor.
Moreover, if the image of an arc in $\DH$ is quasi-radial, then the original arc is quasi-radial
 and the length of the image is smaller than the original length by a definite factor.


\subsection{The saddle fixed point $p_F$ and the proof of Proposition~\ref{p-maximal invariant}}\label{ss-fixed point}
For $ \la \in (0, 1)$ sufficiently close to ~$1$, the  saddle fixed
point  $p_\la = (1 - \la)^{-1} + 1$ of $F_\la$ is contained in the
domain of hyperbolicity $\DH$. So the unstable manifold of $p_\la
$  is a quasi-angular arc while it remains in $\DH.$ As the
unstable cone field is given by $ \cK(z) = \{ z \rho (i +
\varepsilon) \mid \rho \ge 0, \ |\varepsilon| \le 1/3 \},$ we see
that a piece of the unstable manifold of $ p_\la $ is a
quasi-angular arc parameterized by $\theta \mapsto
\rho_\la(\theta) \exp(i \theta)$, with $\rho_\la (0) = p_\la $ and
$\rho_\la (\theta) \ge p_\la \exp(- \theta / 3).$
This remains valid while
$$
p_\la \exp(- \theta / 3)
\ge
\la^{-1} \left( \la - 1 + 4 / \sqrt{1 - \la} \right).
$$
From now on we assume that $ \la \in (0, 1)$ is sufficiently close to ~$1$, so that the function $\rho_\la$ is defined on the interval $[-\pi, \pi]$.

By part~3 of Lemma~\ref{l-unstable cone field}, for every $\theta \in (0, \pi]$
we have $\rho'_\la(\theta) < 0$. By symmetry it follows that for
every $\theta \in [-\pi, 0)$ we have that $\rho_\la'(\theta) > 0$
and that $\rho_\la' (0) = 0$. Therefore the image of $[-\pi,
\pi]$ by the map $\theta \mapsto \rho_\la(\theta) \exp(i\theta)$
forms a Jordan curve, denoted by $\gamma_\la^+$, that is contained
in the unstable manifold of $p_\la$. So $\gamma^+_\la$ is smooth,
except at $\theta = \pm \pi$, which is a point of transversal
self-intersection of the unstable manifold of $ p_\la.$

Note that every map $F$ in $\sFla$ that is sufficiently close to
$F_\la$, has a saddle fixed point $p_F$ that is the continuation
of $p_\la$. Moreover, a piece of the unstable manifold of $p_F$
forms a Jordan curve $\gamma_F^+$ that is smooth, except at a
point of transversal self-intersection. We will denote by $D_F$
the closed disk in $\C$ bounded by the Jordan curve $\gamma^+_F$.

Recall that the maximal invariant set $\Omega_F$ of $F$ is defined
by $\Omega_F = \cap_{m \ge 1} F^m(B_\la)$. In the proof of
Proposition~\ref{p-maximal invariant} below, we will show that
$\Omega_F = D_F \cap F(B_\la^*)$. We say that  $\gamma_F^+$ is the
{\it external boundary of $\Omega_F$}.

\

\proofof{Proposition~\ref{p-maximal invariant}}
In part~$1$ we prove that the the conclusions of the proposition hold for the set $D_F \cap F(B_\la)$ and in part~$2$ we show that $\Omega_F = D_F \cap F(B_\la)$.

\noindent{\bf 1.}
For $F \in \sFla$ sufficiently close to $F_\la$ the set $D_F \cap F(B_\la)$ clearly contains $z = 0$ and $\gamma_F^+$.
So we just need to prove that the set $D_F \cap F(B_\la)$ is homeomorphic to an annulus that is bounded by $\gamma_F^+$ and by a Jordan curve $\gamma_F^-$, and is disjoint from $\gamma_F^-$.

Fix $\varepsilon \in (0, 1)$, put $\gamma^\varepsilon = \{ |z - 1| = 1 - \la + \la \varepsilon \}$ and let $C_F \subset \C$ be the closed annulus bounded by $\gamma_F^+$ and $\gamma^\varepsilon$.
When $F = F_\la$, the set $C_{F_\la}$ is contained in the interior of $F_\la (B_\la)$.
As $F_\la$ is a~$2$ to~$1$ covering between $B_\la^*$ and $F_\la(B_\la)$, it follows that if $F \in \sFla$ is sufficiently close to $F_\la$, then $F$ is a $2$~to~$1$ covering map between $\tC_F = F^{-1}(C_F)$ and $C_F$.
In particular the sets $\tgamma^+_F = F^{-1}(\gamma_F^+)$ and $\tgamma_F^\varepsilon = F^{-1}(\gamma^\varepsilon)$ are Jordan curves.

Denote by $\tD_F^\varepsilon$ the disc bounded by $\tgamma_F^\varepsilon$.
Let $F \in \sFla$ be close to $F_\la$ and put $F = g \circ \tau$, with $g \in \sGla$ close to $\gcan$ and $\tau \in \sTla$ close to $\tau_\la$.
Note that $F(B_\la) = g(\T \times (1 - \la, 2(1 - \la)^{-1} - 1 - \lambda))$.
As the space $\sGla$ is endowed with the strong $C^1$ topology, it follows that for every $F \in \sFla$ sufficiently close to $F_\la$ there is a $C^1$ map $\rho_F^- : \T \to \R$ that is close to the constant map $\theta \mapsto 1 - \la$ and such that
$$
g(\T \times \{ 1 - \la \})
=
\{ \theta \in \T \mid 1 + \rho_F^-(\theta)\exp(4\pi i \theta) \}.
$$
Taking $F$ closer to $F_\la$ if necessary we have
\begin{multline*}
g(\tD_F^\varepsilon)
 =
\{ 1 + t \exp(4\pi i \theta) \mid \theta \in \T, \rho_F^-(\theta) \le t \le 1 - \la + \la \varepsilon \}
= \\ =
\{ 1 + t \exp(2\pi i \psi) \mid \psi \in \T, \min \{ \rho_F^-(\psi/2), \rho_F^-((\psi + 1)/2) \} \le t \le 1 - \lambda + \lambda \varepsilon \}.
\end{multline*}
As $C_F \subset F(B_\la)$, it follows that $D_F \cap F(B_\la)$ is the annulus bounded by $\gamma_F^+$ and by the Jordan curve
$$
\gamma_F^- =  \{ 1 + \min \{ \rho_F^-(\psi/2), \rho_F^-((\psi + 1)/2) \} \exp(2\pi i \psi) \mid \psi \in \T \}.
$$
Clearly $D_F \cap F(B_\la)$ is disjoint form $\gamma_F^-$.

\noindent{\bf 2.}
To prove that $\Omega_F = D_F \cap F(B_\la)$ we first prove that the set $D_F \cap F(B_\la)$ is invariant by $F$ by showing in part~$2.1$ that $F(D_F) = D_F \cap F(B_\la)$ and in part~$2.2$ that $F(D_F \cap F(B_\la)) = F(D_F)$.
This implies that $D_F \cap F(B_\la) \subset \Omega_F$.
Then we complete the proof by showing in part~$2.3$ that $\Omega_F \subset D_F \cap F(B_\la)$.

\noindent{\bf 2.1.}
We keep the notation of part~$1$.
Denote by $\tD_F$ the disc bounded by $\tgamma_F^+$.
As $\tgamma_F^+ = F^{-1}(\gamma_F^+)$ and $F : \tgamma_F^+ \to \gamma_F^+$ is a~$2$ to~$1$ covering map, we have $F(\tD_F) = D_F \cap F(B_\la)$.

We will prove now that $D_F \subset \tD_F$.
Note that $\ell_F^+ = \tgamma_F^+ \cap \gamma_F^+$ is a $C^1$ arc whose image by $F$ is equal to $\gamma_F^+$.
Moreover, the closure of $\tgamma_F^+ \setminus \ell_F^+$ is a $C^1$ arc $\ell_F^-$ having the same end points as $\ell_F^+$, where these arcs intersect transversally.
When $F = F_\la$ we have $\ell_{F_\la}^- = - \ell_{F_\la}^+$, and since $\gamma_\la^+$ is tangent to the cone field $\cK^-$, it follows that $\gamma_F^+ \setminus \ell_{F_\la}^+$ is contained in the interior of $\tD_{F_\la}$.
Now, for $F \in \sFla$ close to $F_\la$ the arc $\ell_F^-$ intersects $\gamma_F^+$ at its $2$ extreme points in a transversal way.
Thus, taking $F$ closer to $F_\la$ if necessary, we have that the arc $\gamma_F^+ \setminus \ell_F^+$ is contained in the interior of $\tD_F$.
Thus $D_F \subset \tD_F$ and $F(D_F) \subset D_F \cap F(B_\la)$.

To prove that $F(D_F) = D_F \cap F(B_\la)$ we just need to prove that $F(D_F) = F(\tD_F)$.
Put $E_F = \tD_F \setminus D_F$ and note that for every $F$ sufficiently close to $F_\la$, we have $E_F \subset \tC_F$.
As $F$ is a~$2$ to~$1$ covering map between $\tC_F$ and $C_F$, to prove that $F(D_F) = F(\tD_F)$ is enough to prove that $F$ is injective on $E_F$.
Let $\tq_F$ and $\tq_F'$ be the end points of $\ell_F^-$.
Note that these points are mapped to the unique point $q_F$ in $\gamma_F^+$ where the unstable manifold of $p_F$ has a transversal self-intersection.
When $F = F_\la$ the set $\ov{E_{F_\la}} \setminus \{ \tq_{F_\la}, \tq_{F_\la}' \}$ is contained in the left half plane, where $F_\la$ is injective.
So, if $F$ is sufficiently close to $F_\la$ it follows that $F$ is locally injective on $\overline{E_F}$.
To prove that $F$ is injective on $E_F$ is enough to prove that a point in $E_F$ near $\tq_F$ and a point in $E_F$ near $\tq_F'$ cannot have the same image.
A piece of the unstable manifold of $p_F$ slightly larger than $\gamma_F^+$, locally separates the plane at $q_F$ into four regions.
The points in $E_F$ near $\tq_F$ and $\tq_F'$ are mapped to opposite regions and are therefore disjoint.
Thus it follows that $F$ is injective on $E_F$ and this completes the proof that $F(D_F) = D_F \cap F(B_\la)$.

\noindent{\bf 2.2.}
We prove now that $F(D_F \cap F(B_\la)) = F(D_F)$.
As for every,
$$
D_\la \setminus F_\la(B_\la^*)
=
\{ z \in \C \mid |z - 1| \le 1 - \la \},
$$
have $-z \in F_\la (B_\la^*)$, it follows that for every $z \in
D_F \setminus F(B_\la^*)$ there is $z' \in D_F \cap F(B_\la^*)$
close to $-z$ such that $F(z') = F(z)$.
Thus $F(D_F) = F(D_F \cap F(B_\la))$.

\noindent{\bf 2.3.}
We now prove that $\Omega_F \subset D_F \cap F(B_\la)$. As $\Omega_F \subset F(B_\la)$ by definition, we just need to prove that $\Omega_F \subset D_F$.
Since $\gamma_F^+$ is a quasi-angular Jordan curve, it follows that exists a constant $C > 0$ such that for every $z \in B_\la \setminus D_F$ the length $L$ of each quasi-radial arc joining $z$ to $\gamma_F^+$ satisfies $\dist(z, D_F) \le L \le C \dist(z, D_F)$.

For a given integer $n \ge 1$, let $z_0, \ldots, z_n \in B_\la \setminus D_F$ be such that for every $j = 0, \ldots, n - 1$ we have $F(z_{j}) = z_{j + 1}$.
Let $\ell_n$ be the radial arc joining $z_n$ to $\gamma_F^+$.
As $F(D_F) \subset D_F$ for $j = 0, \ldots, n - 1$ we can find inductively a quasi-radial arc $\ell_j$ joining $z_j$ to $\gamma_F^+$ and such that $\ell_{j + 1} \subset F(\ell_j)$.
It follows that the length of $\ell_j$ is exponentially small with $j$.
This proves that $F^n(B_\la)$ is contained in an exponentially small neighborhood of $D_F$ and that $\Omega_F \subset D_F$.
\finproof


\section{Fundamental annulus and local unstable manifolds}\label{fundamental annulus}
In this section we construct, for every $\la \in (0, 1)$
sufficiently close to ~$1$ and every $F$ in $\sFla$ sufficiently
close to $F_\la$, an annulus $\AFF$ such that $\AFF \subset F(\AFF)$ (Subsection~\ref{ss-self-covering property}) and such that exists a positive integer $N$ for which $F^N(\AFF) = \Omega_F$ (Subsection ~\ref{ss-covering property}).
The annulus $\AFF$ is defined in
Subsection~\ref{ss-fundamental annulus}. From the property $\AFF
\subset F(\AFF)$ we deduce that for every point $z_0$ in $\AFF$
there is an infinite backward orbit of $F$ in $\AFF$ that starts
at $z_0$ and that each one of these backward orbits  has
associated a local unstable manifold of definite size
(Subsection~\ref{ss-local unstable manifolds}).

\subsection{Fundamental annulus}\label{ss-fundamental annulus}
Fix $\rad \in (0, \exp(-\pi/3))$ and $\la \in (0, 1)$ close to~$1$.
Then put
$$
U_+ = \left\{ z \in \C \mid \Re z \geq 0, \ |z| \geq \rad(1 -
\la)^{-1} \right\},
$$
$$
U_- = \left\{ z \in \C \mid \Re z \leq 0, \ |z - 1| \geq \rad(1 -
\la)^{-1} \right\}
$$
and $ U_0 = U_+ \cup U_-.$ Notice that the boundary $ \gamma_0 $
of $U_0$ is a Jordan curve.
As the Jordan curve $\gamma_\la^+$ is contained in $\{ z \in \C \mid  |z| \ge \exp(-\pi/3) \left( (1 - \la)^{-1} + 1 \right) \}$ (Subsection~\ref{ss-fixed point}), it follows that for every $F$ in $\sFla$ sufficiently close to $F_\la$, the Jordan curve $\gamma_0$ is contained in the disc bounded by $\gamma_F^+$.
Let $ \AFF $ be the annulus bounded by $\gamma^+_F $ and $ \gamma_0 $, that contains $\gamma_F^+$ and is disjoint from $\gamma_0$.
The annulus $ \AFF $ will be called the {\it fundamental annulus} of $F$.

As $F_\la(\gamma_\la^+)$ is contained in $\{ z \in \C \mid |z| \ge \la \exp(-\pi/3)\left( (1 - \la)^{-1} + 1 \right) - \la \}$, it follows that if $\la \in (0, 1)$ is close enough
to ~$1$ and $ F \in \sFla $ is close enough to $ F_\la $, then
\begin{equation}\label{e-internal radious}
F(\gamma_F^+) \subset \{ z \in \C \mid |z| > \rad(1 - \la)^{-1} \},
\end{equation}
and that $\AFF$ and $F(\AFF)$ are both contained in the domain of hyperbolicity $\DH.$

\subsection{Self-covering property of the fundamental
annulus}\label{ss-self-covering property} The purpose of this
subsection is to prove the following proposition.

\begin{proposition}
For $\la \in (0, 1)$ sufficiently close to~$1$ and for $ F \in
\sFla $ sufficiently close to $F_\la$, the annulus $ \AFF $
satisfies $\AFF \subset F(\AFF).$
\end{proposition}

As $\gamma_F^+ \subset F(\gamma_F^+) \subset D_F$ and by \eqref{e-internal radious} the set $F_\la(\gamma_F^+)$ is outside the disc bounded by $\gamma_0$, the proposition is an immediate
consequence of the following lemma.

\begin{lemma}
The closed set $U_0$ is contained in the interior of $F_\la(U_0).$
\end{lemma}

\proof
As
$$
F_{\la}(U_+)
=
\{w \in \C \mid |w - 1| \geq 1 - \la + \la \rad(1 - \la)^{-1} \},
$$
and $1 - \la + \la \rad (1 - \la)^{-1} < \rad (1 - \la)^{-1}$, it follows that $
U_- $ is contained in the interior of $F_{\la}(U_+).$ So it is
enough to show that $U_+$ is contained in the interior of
$F_{\la}(U_-).$

Note that the image by $F_\la$ of the quasi-angular arc $ \{ z \in
\C \mid \Re z \leq 0, |z - 1| = \rad(1 - \la)^{-1} \} $ is a
quasi-angular Jordan curve.
So it is enough to prove that for every $z_0  \in \C$ such that $|z_0 - 1| = \rad(1 - \la)^{-1}$ and $\arg(z_0) \in [\pi/2, 3 \pi / 2]$, we have $|F_{\la}(z_0)| < \rad (1 - \la)^{-1}$.
For such $z_0$ put $ \rho = |z_0| $ and $ \alpha = \arg(z_0)$.
Thus,
$$
\trho = |F_{\la}(z_0) - 1 | = 1 - \la + \la \rho
\ \text{ and } \
\arg(F_{\la}(z_0) - 1) = 2\alpha \mod 2\pi.
$$
We have,
$$
\left( \rad(1 - \la)^{-1} \right)^2
=
1 + \rho^2 - 2\rho \cos\alpha,
$$
$$
|F_{\la}(z_0)|^2
= 1 + \trho^2 - 2\trho \cos(\pi - 2\alpha)
= 1 + \trho^2 + 2\trho \cos 2\alpha.
$$
As $ \rho \ge \rad (1 - \la)^{-1} - 1 > 1, $ we have $ \trho < \rho
$ and $ |F_{\la}(z_0)|^2 < 1 + \rho^2 + 2\rho \cos 2\alpha.$ Hence
$$
|F_{\la}(z_0)|^2
<
\left( \rad(1 - \la)^{-1} \right)^2 + 2\rho (\cos 2\alpha + \cos \alpha).
$$
Since $\alpha \in [\pi/2, 3\pi/2]$ we have $\cos \alpha \le 0$, so
$$
\cos 2 \alpha + \cos \alpha =  (2\cos \alpha - 1)(\cos \alpha + 1) \le 0.
$$
It follows that $|F_{\la}(z_0)| < \rad(1 - \la)^{-1}$, as wanted.
 \finproof

\subsection{Covering property of the fundamental
annulus}\label{ss-covering property} The purpose of this
subsection is to prove the following proposition.

\begin{proposition}\label{p-fundamental annulus}
For every $\la \in (0, 1)$ sufficiently close to~$1$ and every $F$
in $\sFla$ sufficiently close to $F_\la$ there is a positive
integer $ N $ such that $ F^N(\AFF) = \Omega_F$.

\end{proposition}

The proof of this proposition depends on some lemmas.

\begin{lemma}\label{l-scaping orbit} Assume that $\la \in (0, 1)$ is
sufficiently close to~$1$ and that $F$ in $\sFla$ is
sufficiently close to $F_\la$.
Then there is $N \ge 1$ such that for every $z_0 \in F(B_\la^*)$ there exists a positive integer $m_0 \le N$ and a backward orbit $z_0, z_1, \ldots, z_{m_0}$ (i.e. for $j = 1, \ldots, m_0$ we have $F(z_j) = z_{j-1}$), such that $z_{m_0} \in B_\la \setminus F(B_\la^*)$.
\end{lemma}

We will prove this lemma for $ F = F_\la.$ The case when $ F $ is
sufficiently close to $ F_\la$ is an easy consequence.

Let $Q^+ = \{z \in \C \mid \Re z \leq 0 \;{\rm and}\; \Im z \geq 0
\} $ and  $ Q^- = \{z \in \C \mid \Re z \leq 0 \; {\rm and}\; \Im
z \leq 0 \}.$ As $F_{\la} \left( (Q^+ \cup Q^-) \setminus \{ 0 \}
\right) = F_{\la}(\C \setminus \{0\})$, it follows that for every
$ z_0 \in F_{\la}(\C \setminus \{0\})$ there is a pre-image $z_1$
of $ z_0$ by $ F_{\la}$ in $Q^+ \setminus \{0\}$ or in $ Q^-
\setminus \{0\}.$ Moreover, note that $Q^+ \subseteq F_{\la}(Q^-
\setminus \{0\})$ and $Q^- \subseteq F_{\la}(Q^+ \setminus
\{0\}).$ Thus the previous lemma is an immediate consequence of
the following one.

\begin{lemma}
\label{l-preimage growth} Let $z' \in Q^+ \setminus \{0\} $ be
such that $ z = F_{\la}(z') \in Q^-.$ Then the following
properties hold.

\begin{enumerate}
\item[1.] $|z' - 1| \geq \sqrt{2}$. \item[2.] $|z' - 1| - 1 >
\la^{-1} ( |z - 1| - 1).$
\end{enumerate}
\end{lemma}

\proof

\noindent{\bf 1.} If $ |z'| \in (0, 1)$, then $ \Re(F_{\la}(z')) >
0$ so $ F_{\la}(z') \notin Q^-.$ But $ |z' - 1| < \sqrt{2} $
implies $ |z'| < 1 $ \;\,(since $\Re z' \leq 0$), which is not
possible by the previous observation.

\noindent{\bf 2.} Set $|z - 1| = 1 + \mu,$ where $ \mu \geq 0.$
Since $ 1+ \mu = |z - 1| = 1 - \la + \la|z'| $ we have $|z'| = 1 +
\la^{-1} \mu.$ Then observe that, since $\Re z' \leq 0,$ we have
$|z' - 1| > |z'|.$  \finproof

\

\proofof{Proposition~\ref{p-fundamental annulus}}
Let $N$ be the integer given by Lemma~\ref{l-scaping orbit} and for a given $z_0 \in \Omega_F \subset F(B_\la)$ let $m_0 \le N$ and $z_1, \ldots, z_{m_0} \in B_\la$ be as in this lemma.
Thus there is $m \in \{ 1, \ldots, m_0 \}$ such that $z_m \not \in D_F$ and $F(z_m) = z_{m - 1} \in D_F$.
Then~\eqref{e-internal radious} implies that $z_{m - 1} \in \AFF$, and therefore we have $z_0 \in F^{m - 1}(\Omega_F)$.

As $\AFF \subset F(\AFF)$ and $m \ge N$, it follows that $\Omega_F \subset F^N(\AFF)$.
\finproof

\subsection{Local unstable manifolds}\label{ss-local unstable manifolds}
Assume that $\lambda \in (0, 1)$ is close enough to~$1$, so that
the conclusions of Lemma~\ref{l-unstable cone field} hold for every point $z_0$
in $\DH$. It follows that for $F$ sufficiently close to
$F_\lambda$, the set
$$
\Gamma_F = \left\{ \{ z_i \}_{i \ge 0} \mid z_i \in \DH \text{ and
} F(z_{i + 1}) = z_i, \text{ for } i \ge 0 \right\}
$$
is non empty and that for every $ z_0 \in \AFF $ there is an
infinite backward orbit  in $\Gamma_F $ starting from $ z_0.$ From
the theory of hyperbolic sets we know  that there exists $\alpha >
0$ such that for every infinite backward orbit  $ \underline{z} $
in $\Gamma_F $ starting from a point $ z_0 \in \AFF,$  the set
\begin{multline*}
W^u_\alpha(\underline{z}) = \left\{ w_0 \in \DH \right. \mid
w_0 \text{ has a backward orbit $\underline{w} \in \Gamma_F$} \\
\left. \text{with $|w_m - z_m| < \alpha$ for all $m \ge 0$}
\right\}
\end{multline*}
is a quasi-angular arc. The set $W^u_\alpha (\underline{z})$ will
be called the {\it local unstable manifold of $ \underline{z}$}.
Moreover, for every $\varepsilon > 0$ there exists $ \delta > 0 $
and a positive integer $ N $ such that, if $\underline{z}'$ is an
infinite backward orbit in $\Gamma_F$ such that for every $0 \le m
\le N$ we have $| z'_m - z_m| < \delta$, then the $C^1$ distance
between the local unstable manifolds $W^u_\alpha(\underline{z})$
and $W^u_\alpha(\underline{z}')$ is at most $\varepsilon$.

\section{Topological dynamics on the attractor}\label{dinamica atractor}
In this section we prove Theorem~\ref{t-topological dynamics}
(Subsection~\ref{proof-top dynamics}), by showing that for every
$\la \in (0, 1)$ sufficiently close to~$1$ and every area expanding map $F$ in
$\sFla$ sufficiently close to $F_\la$, the stable manifold of
the saddle fixed point $p_F$ is dense in $B_\la$
(Subsection~\ref{ss-stable density}) and that the unstable
manifold of $p_F$ is dense in  $\Omega_F$
(Subsection~\ref{ss-unstable density}).

Throughout this section we fix $\la \in (0, 1)$ close to ~$1$ and
an area expanding map $F$ in $\sFla$ close to $F_\la$.

\subsection{Density of the stable manifold of $p_F$}\label{ss-stable density}
Consider a connected  open subset $V$ of $B_\lambda$. Suppose for a moment
that there is a positive integer $m$ such that for every $k = 1, \ldots, m$ the map $F$ is injective on $F^k(V)$.
As $F$ is area expanding, it follows that the area of $F^k(V)$ grows exponentially fast with $k$, as long as $k \le m$.
We conclude that there is an integer $m$ for which $F$ is not injective on $F^m(V)$. The following lemma
implies that in this case $F^m(V)$, and hence $V$, intersects  the
stable manifold of $p_F$. From this it follows that the stable
manifold of $p_F$ is dense in $B_\lambda$.

\begin{lemma}\label{l-stable density}
Every connected open subset of $B_\lambda$ where $F$ is not injective
intersects the stable manifold of $p_F$.
 \end{lemma}

\proof
Recall that $F_\la = \gcan \circ \tau_\lambda$, where $\tau_\lambda \in \sTla$ and $\gcan \in \sGla$ are defined in Subsection~\ref{ss-the spaces}.
Note that $F_\lambda$ maps the interval $(0, 2(1 - \lambda)^{-1} + 1) \subset \R \subset \C$ onto $(2 - \lambda, 2(1 - \la)^{-1})$.
It follows that $(0, 2(1 - \lambda)^{-1} + 1)$ is contained in the stable
manifold of $p_\lambda$.
On the other hand, note that the preimage of the arc $(1, 2(1 -
\lambda)^{-1})$ by $\gcan$ consists of the arcs $\{ 0 \} \times
I$ and $\{ 1/2 \} \times I$, where $I \subset \R$ is the interval
$I = (0, 2(1 - \la)^{-1} - 1)$, and that each one of these arcs
intersects transversally $\T \times \{ 1 - \la \}$ and
$\T \times \{ 2 (1 - \lambda)^{-1} - 1 - \lambda \}$.
Moreover these arcs separate $\T \times [1 - \lambda, 2(1 - \lambda) - 1 - \lambda]$ into~2 sets where $g_*$ is injective.

Given $F \in \sFla$ close to $F_\lambda$, let $\tau \in \sTla$ be close to $\tau_\la$ and $g \in \sGla$ be close to $\gcan$, such that $F = g \circ \tau$.
Denote by $\gamma^s$ an arc of the stable manifold of $p_F$ that is $C^1$ close to $(1, 2(1 - \lambda)^{-1})$.
Since $g$ is $C^1$ close to $\gcan$ in the $C^1$ topology, the preimage of $\gamma^s$ by $g$ contains~2 arcs $\gamma$ and $\gamma'$, that intersect transversally the circles $\T \times \{ 1 - \la \}$
and $\T \times \{ 2 (1 - \lambda)^{-1} - 1 - \lambda \}$. Moreover,
these arcs separate the cylinder $\T \times [1 - \lambda, 2(1 - \lambda)^{-1} - 1 - \lambda]$
into two parts, on each of which $g$ is injective.

It follows that $\tau^{-1}(\gamma \cup \gamma') = F^{-1}(\gamma^s)$ divides $B_\lambda^*$ into two parts, on each of which $F$ is injective.
So every open and connected subset of $B_\lambda$ where
$F_\lambda$ is not injective must intersect $F^{-1}(\gamma^s)$. As this last set is contained in the stable manifold of $p_F$, the statement follows.
\finproof

\subsection{Density of the unstable manifold of $p_F$}\label{ss-unstable density}
Together with Proposition~\ref{p-fundamental annulus}, the
following lemma implies that the unstable manifold of $p_F$ is
dense in $\Omega_F$. In fact it is easy to see that the
local unstable manifold of each infinite backward orbit in
$\Gamma_F$ is approximated in the $C^1$ topology by arcs of the
unstable manifold of~$p_F$.

\begin{lemma}\label{radial arc}
Every quasi-radial arc contained in $ \AFF $ intersects
transversally the unstable manifold of $ p_F$ at some point.
\end{lemma}

\proof Given a quasi-radial arc $\gamma_0$ contained in $\AFF$, choose $z_0 \in \gamma_0$ and let $\underline{z} = \{ z_m \}_{m \ge 0} \in \Gamma_F$ be an infinite backward orbit of points in $\AFF$, starting at $z_0$.
For $j \ge 1$ define inductively a quasi-radial arc $\gamma_j$ joining $z_j$ to some point in $\gamma_F^+$ and such that $F(\ell_j)$ is contained in $\ell_{j - 1}$.
As for every $j \ge 0$ the point $z_j$ belongs to $\AFF$, the length of $\ell_j$ is bounded independently of $j$.
It follows that the length of the arc $F^j(\ell_j)$ is exponentially small with $j$.
Thus $ F^j(\gamma_j \cap \gamma^+_F)$ is a point of the unstable manifold of $p_F$ that is exponentially close to $z_0$.
\finproof


\subsection{Proof of  Theorem~\ref{t-topological dynamics}}\label{proof-top dynamics}
Let $ U, V $ be open subsets of $B_\la$ that intersect the
attractor set $ \Omega_F$ of $F$. We now prove that there exists a
positive integer $ m_0 $ such that for every $m \ge m_0$ we have
$F^m(U) \cap V \neq \emptyset.$

As the unstable manifold of $p_F$ is dense in $\Omega_F$, there
exists a small open set $V_0$ and a positive integer $j_0$ such
that $V_0$ intersects the local unstable manifold of $p_F$ and
such that $F^{j_0}(V_0) \subset V.$ As the stable manifold of
$p_F$ is dense in $\Omega_F$, there exists a positive integer
$k_0$ such that $ F^{k_0}(U) $ intersects the local stable
manifold of $p_F$. By the Inclination Lemma, there exists a
positive integer $\ell_0$ such that for every $\ell \ge \ell_0$ we
have $F^\ell(F^{k_0}(U)) \cap V_0 \neq \emptyset$.

Taking $ m_0 = j_0 + \ell_0 + k_0 $ we see that for every $m \ge
m_0$ we have $F^m(U) \cap V \neq \emptyset$. This proves
Theorem~\ref{t-topological dynamics}. \finproof


\section{Eventually onto property and periodic points}\label{s-periodic points}
This section is devoted to the proof that the set of periodic sources and the set of periodic saddles of a generic area expanding map $F$ in $\sFla$ that is close to $F_\la$, are both dense in~$\Omega_F$ (Proposition~\ref{p-periodic points} in Subsection~\ref{ss-attractor set}).
The proof is based on the following intermediate property.
\begin{definition}
  \label{d-eop}
  We say that $F \in \sFla$ has the {\sf eventually onto property} if for every open subset $V$ of $\Omega_F$
$$
\interior(\Omega_F) \subset \cup_{m \ge 0} F^m(V).
$$
\end{definition}
We first prove that the conclusions of Proposition~\ref{p-periodic points} hold for every map $F$ close to $F_\la$ with the eventually onto property and satisfying a certain genericity condition (Proposition~\ref{p-sources} in Subsection~\ref{ss-sources}).
This genericity condition is expressed in terms of the eigenvalues of the derivative of $F$ at a fixed source $p^+_F$ of $F$ (see Subsection~\ref{ss-p+} below).

We then show that
for every $ F $ in $\sFla$ near $F_\lambda$, the eventually onto property
holds for all open sets $V$ containing $p^+_F$
(Proposition~\ref{p-eop at p+} in Subsection~\ref{ss-eop at p+}) and we complete the proof of Proposition~\ref{p-periodic points} in Subsection~\ref{ss-iterated preimages of p+}, by showing that the eventually onto property holds for a generic area expanding map near~$F_\la$.

\subsection{Limit behavior and the fixed source $p^+_F$}\label{ss-p+}
When $\lambda \rightarrow 1$ the map $F_\lambda$ converges in the
$C^\infty$ topology to the map $\Gcan: \C^* \to \C$ defined by,
$$
\Gcan(z) = |z| (z/|z|)^2 + 1.
$$
The map $\Gcan$ extends continuously to $z=0$ by setting $\Gcan(0) = 1$.
Moreover, note that $\Gcan$ has constant Jacobian equal to $2$, and
that $\Gcan$ is injective on the upper half plane $\H = \{ \Im z > 0
\}$.

It is easy to see that  the point $p^+ = \exp(\pi i/3)$ is the unique
fixed point of $\Gcan$ in the upper half plane. Moreover, the
derivative of $\Gcan$ at $p^+$ has complex eigenvalues (not in $\R$) of modulus
larger than~$1$, so that $p^+$ is a hyperbolic source.

It follows that for $\lambda \in (0, 1)$ sufficiently close to~$1$ and $F$
 sufficiently close to $F_\lambda$, the map $F$ has a (unique) fixed source $p^+_F$ near $p^+$.
Moreover the derivative of $F$ at $p^+_F$ has complex eigenvalues (not in
$\R$).

\subsection{Dynamics of maps with the eventually onto property}\label{ss-sources}
This subsection is devoted to the proof of the following proposition.
 Recall that the homoclinic class of a saddle periodic point is, by definition, the closure of the set formed by the transversal intersections between the stable and unstable manifolds of the saddle periodic point.
Every point in a homoclinic class is accumulated by saddle periodic points, see Remark~\ref{r-saddles}.
\begin{proposition}\label{p-sources}
 There exists a neighborhood $\sU$ of $F_\la$ in $\sFla$ such that for every $F \in \sU$ with the eventually onto
property the following hold:

\begin{enumerate}
\item[1.] The periodic sources of $F$ are dense in $\Omega_F$.

\item[2.] If the arguments of the eigenvalues of $p^+_F$ are
irrational multiples of $\pi$, then the homoclinic class of $p_F$
is equal to $\Omega_F$. In particular, the periodic saddles of $F$ are
dense in $\Omega_F$.
\end{enumerate}
\end{proposition}
\proof
Let $\la \in (0, 1)$ be close to~$1$ and let $F$ be a map in
$\sFla$ close to $F_\la$ and that satisfies the eventually onto
property. In particular, the iterated preimages of $p_F^+$ are
dense in $\Omega_F$.

Let $ V $ be a small open neighborhood of $ p^+_F $ such that $ F
$ has a local inverse $f$ defined on $V$, such that $f(\ov{V}) \subset V$ and such that $f$ is uniformly contracting, so that $ \cap_{m \ge 0} f^m(V) =
\{p^+_F\}. $

\noindent{\bf 1.} Let us prove that the set of periodic sources of
$F$ is dense in $\Omega_F.$ Let $ q \in \Omega_F $ and $k \ge 1$
be such that $ F^k(q) = p^+_F.$ It is enough to prove that close
to $ q$ there is a periodic source of $ F.$

Let $ U $ be a small open neighborhood of $ q. $  Let $ \{ q_m
\}_{m \ge 1}$ be such that for every $m$ we have that $F^m(q_m) =
q$ and such that $q_m \to p^+_F $ as $ m \to \infty $. Let $ m_0$
be such that for every $m \ge m_0$ we have $q_m \in V.$  Let $
U_0$ be a small open neighborhood of $ q_{m_0} $ contained in $V$,
such that $F^{m_0}(U_0) \subset U$, $F^{m_0 + k}(U_0) \subset V$
and such that for every $j = 1, \ldots, m_0 + k$ we have that
$F^j: U_0 \to F^j(U_0)$ is a diffeomorphism.

Notice that $ F^{m_0 + k}(U_0)$ is an open neighborhood of $
p^+_F.$ Finally, take $ \ell \ge 1$ large enough so that $
\overline{f^\ell(U_0)} \subset F^{m_0 + k}(U_0)$ and so that
$F^{\ell + m_0 + k} $ restricted to $ W = f^\ell(U_0) $ is
uniformly expanding.  It follows that $ \overline W \subset
F^{\ell + m_0 + k}(W) $ and that $ W$ contains a periodic source $
\widetilde{q}$ of $F$ of period $ \ell + m_0 + k $.
Thus, $F^{\ell + m_0}(\widetilde{q}) \in U$ is a periodic source for $ F
$ close to $ q.$

\noindent{\bf 2.} Assume that the argument of the eigenvalues of
$D_{p_F^+}F$ are irrational multiples of $\pi$.

As the stable manifold of the saddle fixed point $p_F$ of $F$ is
dense in $\Omega_F$ it follows that there is a point $z_0$ in $V$
contained in the stable manifold of $p_F$. Let $v_0 \in T_{z_0}
B_\la$ be a vector that is tangent to the stable manifold of $p_F$
at $z_0$. For $m \ge 1$ put $z_m = f^m(z_0)$ and $v_m = D_{z_0}f^m
v_0 \in T_{z_m} B_\la$. As the argument of the eigenvalues of
$D_{p_F^+} F$ are irrational multiples of $\pi$, it follows that
the arguments of the $v_m$ are dense in $[-\pi, \pi]$.

So, each iterated preimage $q$ of $p_F^+$ in $A_F$ is accumulated
by arcs of the stable manifold that are quasi-radial. Hence
Lemma~\ref{radial arc} implies that $q$ is accumulated by points
of transversal intersection between the stable and unstable
manifold of $p_F$. From Proposition ~\ref{p-fundamental annulus}
it follows that the transversal homoclinic intersections of $p_F$
are dense in $\Omega_F.$ \finproof

\begin{remark}\label{r-saddles}
{
We now prove that arbitrarily close to a transversal
homoclinic point of $p_F$ there is a saddle periodic point of $F$.
Although this fact is well-known in the case of diffeomorphisms, some care should be taken for the singular endomorphisms considered here.

Let $ U $ be an open subset of $\Omega_F$ and let $ q \in U $ be a
point of transversal intersection between the stable and unstable
manifolds of $p_F.$ Then, for every positive integer $ m $ we have
$ q_m = F^m(q) \in \C^* $  and $ q_m \rightarrow p_F $ as $ m
\rightarrow \infty.$ Consider an infinite backward orbit $
\{\overline{q}_m\}_{m \geq 0} \subset \C^* $ of $ q $ (that is $
\overline{q}_0 = q $ and $ F(\overline{q}_{m+1}) = \overline{q}_m
$ for every $ m $) such that $ \overline{q}_m \rightarrow p_F $ as
$ m \rightarrow \infty.$

Let $ N $ be a large integer, such that $ q_N $ (resp. $
\overline{q}_N $) is in the local stable (resp. unstable) manifold
of $ p_F.$ There is a small open neighborhood $ \overline{V}_N $
of $ \overline{q}_N$ whose boundary are two small quasi-radial
arcs and two quasi-angular arcs, such that the following
properties hold.

\begin{enumerate}

\item[(i)] For every $ m, \, 0 \leq m \leq 2N, F^m(\overline{V}_N)
\subset \C^* $  and $ V_0 = F^N(\overline{V}_N) $ is a small
neighborhood of $ q$ that has the shape of a rectangle with two
sides `parallel' to the piece of stable manifold of $p_F$ that
contains $q,$ and two sides `parallel' to the piece of the
unstable manifold of $p_F$ that contains $ q$,

\item[(ii)] $V_N = F^{2N}(\overline{V}_N) = F^N(V_0)$ is a small
neighborhood of $q_N$ that has the shape of a rectangle with two
sides 'parallel' to the local stable manifold of $p_F$, and two
sides that are transversal to the local stable manifold of $p_F.$

\end{enumerate}

Clearly, the two sides of $ \overline{V}_N $ that are quasi-radial
arcs map onto the two sides of $ V_N $ that are parallel to the
local stable manifold of $ p_F.$ We now take a sufficiently large
integer $M$  so that $ F^M(V_N) \cap \overline{V}_N $ is a thin
strip that crosses from one quasi-radial arc in the boundary of $
\overline{V}_N $ to the other.

It follows that $ F^M(V_N) \cap \overline{V}_N \ne \emptyset $ and
therefore $V_0$ contains a saddle periodic point of $ F.$}
\end{remark}

\subsection{Eventually onto property at $p^+_F$}\label{ss-eop at p+}
This subsection is devoted to prove the following proposition.

\begin{proposition} \label{p-eop at p+}
For all  $\lambda \in (0, 1)$ sufficiently close to $1$ and
every $F$ sufficiently close to $F_\lambda \in \sFla$, the following
property holds. For every  neighborhood $U$ of the fixed point
$p^+_F$,
$$
\interior (\Omega_F) \subset \cup_{m \ge 0} F^m(U).
$$
\end{proposition}

The following lemma reduces this proposition to show that a
certain arc $I$ in $\C$ is contained in the basin of the fixed
source $p^+_F$. This last fact is proven in lemmas
~\ref{preliminar disk} and ~\ref{basin of p^+}.

\begin{lemma} Let $V$ be a neighborhood of the arc $I = \{ t \in \C \mid t \in \R, t \in [0, 1] \}$.
If $\lambda \in (0, 1)$ is sufficiently close to~$1$ and $F$ is
sufficiently close to $F_\lambda$, then $\interior(\Omega_F)
\subset \cup_{m \ge 0} F^m(V)$.
\end{lemma}

\proof

It is enough to prove that
$$
\interior(\AFF) = \AFF \setminus \gamma^+_F \subset \cup_{m \ge 0} F^m(V).
$$
In fact, by Proposition~\ref{p-fundamental annulus} it
follows that there exists an integer $N \ge 1$ such that
$F^N(\interior(\AFF)) = \interior (\Omega_F)$.

After some preliminary considerations in part~$1$ this is proven in part~$2$.

\partn{1}
First notice that if  $ \la \in (0, 1)$ is close enough to $ 1$,
then $F_\la(V) \cup F_\la^2(V)$  contains the fundamental domain
$[2, 2 + \la] \times \{0\}$ of the stable manifold of $p_\lambda$.
Therefore, $\cup_{m=0}^{m = 22} F_\la^m(V) $ contains a
neighborhood $V_0$ of a fundamental domain $D_0 \subset \R$ of the stable
manifold of $p_\lambda$ which is contained in $\{ z \in \C \mid |z| \ge 20 \}$.
For every $F$ close to $F_\la$ the open set $V_0$ is also a
neighborhood of a fundamental domain $D_0(F)$ of the stable
manifold of $p_F.$

Recall that if $\la $ is sufficiently close to ~1, then the
unstable cone  fields of $F_\la$ (see Subsection~\ref{ss-unstable cone
field}) are defined and invariant  in $\{ z \in \C \mid |z| \ge 20 \}$ (Lemma~\ref{l-unstable cone field}). The same happens for every $F $
close to $F_\la.$ In particular, if $x$ is a real positive number
satisfying $x \ge 20 \exp({\pi /3})$, then a quasi-angular arc through $x$
touches the negative real axis before intersecting the circle $\{ z \in \C \mid |z| = 20 \}$.
We will assume that if $x \in D_0$ then $ 20 \exp({\pi /3}) \le x \le
180.$

\partn{2}
At a point $x_0 \in D_0$ let $\eta$ be a quasi-angular arc of
length $\ell$ that contains $x_0$.
Then, since for every positive real number $x > 0$ we have, $F_\la(x) < x + 2,$ it follows by induction that for every positive integer $m$ we have $F_\la^m(x_0) < x_0 + 2m < 180 + 2m $ and $F_\la^m(\eta) $
is quasi-angular arc of length $2^m \ell$  or contains a
quasi-angular Jordan curve around the origin.

For  $m_0$ so that $2^{m_0} \ell > 180 + 2m_0  $ we have
 that $F_\la^{m_0}(\eta) $ contains  a quasi-angular Jordan
curve around the origin.

Now take a small rectangle shaped region $W_0$ contained in $V_0$
and foliated by quasi-angular arcs trough points in $D_0$ (so that
$D_0$ crosses $W_0$ from one side to the opposite one) and such
that the image of the arc through the endpoint of
$D_0$ closer to the origin contains the arc through the other endpoint of $D_0$.

From the above argument, there exists a positive integer $m_0$
such that $F_\la^{m_0}(W_0)$ contains an annulus $A_{m_0}$ around
the origin. Notice that the image by $F_\la $ of the portion of the
internal boundary of $A_{m_0}$ with positive real part is exactly
the external boundary of $A_{m_0}$. This $m_0$ is the same for
every $\la $ close enough to ~1. And by continuity a similar
construction is valid for $F $ close enough to $F_\la.$

Choose $\la $ sufficiently close to ~$1$ in such a way that
  the fundamental annulus $\AFF$
is disjoint from  $A_{m_0} $ and contained in the unbounded
component of the complement of $A_{m_0} $.
For every $m \ge m_0$ inductively define  an annulus $A_m$~by
$$A_{m+1} = F_\la(A_m \cap \{z \in \C \mid \Re z \ge 0\}).$$
Similarly define the annulus $A_m(F), m \ge m_0,$ for every $F $
close enough to~$F_\la.$

By the Inclination Lemma it follows that $$ \AFF
\setminus \gamma_F^+ \subset \cup_{m \ge m_0} F(A_m(F)) \subset
\cup_{m \ge 0}F^m(V_0).$$ \finproof

\begin{lemma}\label{preliminar disk} There is an open subset $D$ of the upper half plane $\H$,
that is bounded
 by a Jordan curve and such that $\Gcan(D)$ contains both the closure of $D$ and the
 arc~$I.$
\end{lemma}

\proof Given $\varepsilon > 0$ small, set
$$
D' = \left\{ z \in \C^* \mid   \arg(z) \in (\varepsilon, \pi / 2),
\,\, |z - p^+| < |2i - p^+| \right\}.
$$
\noindent {\bf 1.} We will prove that if $\varepsilon > 0$ is
sufficiently small, then $\ov{D'} \cap \Gcan(\partial D') = \{ 0 \}$.
In part~2 below we conclude the proof of the lemma from this fact.

The boundary of $D'$ consists of the arc $J_1 = \{ it \mid t \in
\R, 0 \le t \le 2 \}$, an
 arc of the form $J_2 = \{ \rho \exp (i \varepsilon) \mid 0 \le \rho \le \rho_0 \}$,
 for some $\rho_0 > 0$, and an arc $J_3$ of the
  circle $\{ z \in \C \mid |z - p^+| = |2i - p^+| \}$.

Clearly, $\Gcan(J_1) = \{ t \in \C \mid t \in \R, -1 \le t \le 1 \}$
intersects $\ov{D'}$ only at $0$. Let us show that $\Gcan(J_2) = \{
\rho \exp( 2i \varepsilon) + 1 \mid 0 \le \rho \le \rho_0 \}$
 is disjoint from $\ov{D'}$.
As $\arg(1 + \exp (2i\varepsilon)) = \varepsilon$, it follows that
for $0 \le \rho < 1$ the point $1 + \rho \exp(2i\varepsilon)$ is
disjoint from $\ov{D'}$. As $|\, 2i - p^+| < |\, 2 - p^+| $, it
follows that if $\varepsilon$ is sufficiently small, then for all
$\rho \ge 1$
$$
|\, 1 + \rho \exp(2i\varepsilon) - p^+| \ge |\, 1 + \exp(2i
\varepsilon) - p^+| > |\, 2i - p^+|,
$$
so that $\rho \exp(2i\varepsilon) + 1 \not \in \ov{D'}$.

To see that $\Gcan(J_3)$ is disjoint from $\ov{D'}$, observe that for
each pair of
 points $z, z' \in \C^*$ such that $ \arg(z), \arg(z') \in [0, \pi/2]$, we have
  $|\Gcan(z) - \Gcan(z')| \ge |z - z'|$, with equality if and only if $z/z' \in \R$.
So, for every $z \in J_3$ we have $|\Gcan(z) - p^+| \ge |z - p^+|$,
with equality if and only if  $z = (1 + |\, 2i - p^+|)\exp(\pi i
/3)$. Is easy to check that the image of this point has negative
real part. We conclude that $\Gcan(J_3)$ is disjoint from $\ov{D'}$.

\noindent {\bf 2.} For $\delta > 0$ small set $D'' = \{ z \in D'
\mid |z| > \delta \}$. Note that $D''$ is bounded by a Jordan
curve and $\ov{D''} \subset \H$. We will show that, if $\delta >
0$ is sufficiently small, then $\ov{D''} \cup I \subset \Gcan(D'')$.
Then any sufficiently small neighborhood $D$ of $\ov{D''}$ that is
bounded by a Jordan curve, will satisfy the desired properties.

As $\Gcan(0) = 1$ does not belong to $\ov{D'}$, if $\delta > 0$ is
sufficiently small, then $\ov{D''} \subset \Gcan(D'')$. When $\delta <
1$, the arc $I$ is clearly contained in the closure of $\Gcan(D'')$.
\finproof

\

By the previous lemma it follows that, if we take $\lambda$
sufficiently close to~$1$, then every map $F$ sufficiently close
to $F_\lambda$ satisfies the following properties.
\begin{enumerate}
\item[1.] $F$ is injective on $\ov{D}$ and $\ov{D} \subset F(D)$.
\item[2.] The Jacobian of $F$ on $D$ is larger than a constant
larger than one. \item[3.] The fixed source $p^+_F$ is contained
in $D$ and it is the unique fixed point of $F$ in $\ov{D}$.
Moreover, the derivative of $F$ at $p^+_F$ has complex eigenvalues  (not in
$\R$) of norm larger than~$1$.
\end{enumerate}
\

\begin{lemma}\label{basin of p^+} Let $\lambda \in (0, 1)$ be close to~$1$ and let $F$ be an endomorphism in $\sFla$  close to
$F_\lambda$ satisfying the properties above. Denote by $f:F(D) \to
D$ the inverse of $F$ restricted to $F(D)$. Then
$$
\cap_{m \ge 1} f^m(D) = \{ p_F^+ \}.
$$
In particular, for every neighborhood $U$ of $p^+_F$ there is an
integer $m$ such that $D \subset F^m(U)$.
\end{lemma}

 \proof Set $K =
\cap_{m \ge 1} f^m(D)$. As the Jacobian of $f = F^{-1}|_{F(D)}$ is
smaller than a constant smaller than~$1$, it follows that $K$ has
~0 Lebesgue measure, and hence empty interior. Assume by
contradiction that $K$ is not equal to $\{ p^+_F \}$. We will
prove then that $f$ has a fixed point in $D$ distinct from
$p^+_F$. This contradicts property~$3$ above and proves the lemma.

Let $z_0 \in K$ be different from $p^+_F$ and denote by $K_0$ the
set of accumulation points of the forward orbit of $z_0$ by $F$.
So $K_0 \subset K$ and $F(K_0) = K_0$. Moreover, as $p^+_F$ is a
source, $p^+_F \not \in K_0$.

Let $h : \D \to D$ be a linearizing coordinate of $F$ near
$p^+_F$, so that $h(0) = p^+_F$ and so that for every $w \in \D$
satisfying $|w| < 1/2$ we have $F(h(w)) = h(2w)$. The forward
orbit by $f$ of each point in $h(\D)$ converges to $p^+_F$. As
$f(K_0) = K_0$ and $p_F^+ \not \in K_0$, it follows that $K_0 \cap
h(\D) = \emptyset$.

Since $K$ has empty interior, there is $w_0 \in \D$ such that
$h(w_0) \not \in K$. Let $m \ge 1$ be such that $h(w_0) \not \in
f^m(D)$ and set
$$
t_0 = \inf \left\{ t \mid h(tw_0) \not \in f^m(D) \right\} \
\text{ and } \ \gamma = h \left( \left\{ tw_0 \mid  t \in [0, t_0]
\right\} \right).
$$
As $D$, and hence $f^m(D)$, are bounded by a Jordan curve, it
follows that $D_0 = f^m(D) \setminus \gamma$ is homeomorphic to a
disk. Moreover, as $f(\gamma) \subset \gamma$, we have $f(D_0)
\subset D_0$, and since $K_0$ is disjoint from $\gamma \subset
h(\D)$, we have $K_0 \subset D_0$. Let $\tf : D_0 \to D_0$ be a
homeomorphism that coincides with $f$ on $f(D_0)$, such that
$\tf(\gamma) = \gamma$ and such that every point in $D_0$ enters
$f(D_0)$ under forward iteration by $\tf$. In particular $\tf$
does not have fixed points. But, since $K_0$ is a compact subset
of $D_0$ that is invariant by $\tf$, Brouwer's translation theorem
applied to $\widetilde{f}$ implies that $\tf$ has a fixed point in
$f(D_0)$. So we get a contradiction that proves the lemma.
\finproof

\subsection{Proof of Proposition~\ref{p-periodic points}}\label{ss-iterated preimages of p+}
In view of Proposition~\ref{p-eop at p+}, a map $F$ in
$\sFla$ close to $F_\la$ satisfies the eventually onto property
if and only if the iterated preimages of $p_F^+$ are dense in
$\Omega_F$.
So Proposition~\ref{p-periodic points} is a direct consequence of Proposition~\ref{p-sources} and the following proposition.

Recall that for an open set $\sU$ of $\sGla$ and $\tau \in \sTla$ we denote
$$
\sU_\tau = \{ g \in \sU \mid F = g \circ \tau \text{ is area expanding}\}.
$$
\begin{proposition}\label{p-eop}
There is a neighborhood $\sU$ of $\gcan$ in $\sGla$ such that for every $\tau \in \sTla$ there is a residual subset $\sR'_\tau$ of
$\sU_\tau$ with the following property. For every $g \in \sU_\tau
$ the map $F = g \circ \tau$ is such that the iterated preimages
of the fixed source $p_F^+$ are dense in $B_\la$. In particular,
$F = g \circ \tau$ satisfies the eventually onto property.
\end{proposition}
\proof
Let $\sU$ be a sufficiently small neighborhood of $\gcan$ in $\sGla$ such that for every $\tau \in \sTla$ sufficiently close to $\tau_\lambda$, all the results of Section~\ref{technical results} hold for $F = g \circ \tau$.

Given an open subset $U$ of $B_\la$, put
$$
\sA(U) = \left\{ g \in \sU_\tau \mid
p^+_F \in \cup_{m \geq 0}F^m (U)  \text{ where } F = g \circ \tau  \right\}.
$$
Clearly $\sA(U)$ is an open subset of $\sU_\tau$. We will show
that $ \sA(U) $ is dense in $\sU_\tau.$ The residual set
$\sR'_\tau$ will be the intersection of the sets $\sA(U)$, where
$U$ runs through a countable basis for  the topology of $B_\la$.

Given $g \in \sU_\tau$ put $F = g \circ \tau$ and consider an
iterated preimage $ q_0 $ of $p^+_{F}$ in  the fundamental annulus
$ \AFF. $ Consider a $C^1 $ family of functions $\{ h_\varepsilon
\}$ in $\sU_\tau$, such that $h_0 = g$, such that $h_\varepsilon$
coincides with $ g $ on $ \{ |z| \geq 15 \} $ and such that, if we
denote by $ q_\varepsilon $ the iterated preimage of
$p_{F_\varepsilon}^+$ by $F_\varepsilon = h_\varepsilon \circ
\tau$ that is the continuation of $q$, then $q_\varepsilon$ moves
following a radial arc.

Note that for every small $\varepsilon > 0$ the maps
$F_\varepsilon$ and $F_0$ coincide on $\DH$ and we have
$p_{F_\varepsilon} = p_{F_0}$. As the stable manifold of the
saddle fixed point $p_{F_0}$ is dense in $B_\la$, there exists $ m
$ such that $ F_0^m (U) $ intersects the local stable manifold of
$p_{F_0}$.

Let $ \ell \subset \AF_{F_0}$ be a smooth arc contained in $F_0^m
(U)$ that intersects transversally the local stable manifold of $
p_{F_0}$. Shrinking $\ell$ if necessary, we assume that for every
sufficiently small $\varepsilon > 0$ we have $ \ell \subset
F^m_\varepsilon(U)$. From the Inclination Lemma and from Lemma
~\ref{radial arc}, it follows that  there exists a sub-arc $
\ell'$ of $\ell $ and a positive integer $ N $, such that for
every $j = 0, \ldots, N$ the arc $ F_0^j(\ell') $ is a quasi-angular
arc contained in $ \DH $ and such that $F_0^N(\ell')$ intersects
transversally, and in a non empty way, the radial arc defined by
$q_\varepsilon.$ So, for some $\varepsilon \in (-\varepsilon_0,
\varepsilon_0)$, we have $q_\varepsilon \in F_0^N(\ell')$. As for
every $j = 0, \ldots, N$ the arc $F_0^j(\ell')$ is contained in
$\DH$, where $F_0$ coincides with $F_\varepsilon$, we have
$F^j_\varepsilon(\ell') = F_0^j(\ell')$. Thus $q_\varepsilon \in
F_\varepsilon^N(\ell') \subset F_\varepsilon^{N + m}(U)$.
\finproof

\section{Robust tangencies and wild hyperbolic sets}\label{robust tangencies}
In this section we prove Theorem~\ref{t-wild for endos}.
Part~$1$ of this theorem is a direct consequence of Lemma~\ref{wild} in Subsection~\ref{wildhypset}.
The proof of part~$2$ is based on a strong property: every ``curved'' arc near $\AFF$ is tangent to the unstable manifold of some infinite backward orbit in
$\Gamma_F$ (Proposition~\ref{curved arcs and tangencies}).
We complete the proof of part~$2$ of Theorem~\ref{t-wild for endos} in Subsection~\ref{the stable manifold is curved} by showing that there is an arc of the stable manifold of the saddle fixed point $p_F$ of $F$ that is ``curved''.

Throughout this section we assume that $\lambda \in (0, 1)$ is
close enough to~$1$, so that all the properties of Section~\ref{technical results} are satisfied.

\subsection{Wild hyperbolic set}\label{wildhypset}
Fix $F$ in $\sFla$ close to $F_\la$.
Then the set
$$
 W_F = \{z \in \DH \mid F^m(z) \in \DH \;
\; \text{ for every } \; m \ge 1\},
$$
is an uniformly hyperbolic forward invariant set for~$F$.

The saddle fixed point $p_\la$ of~$F_\la$ is contained in $W_{F_\la}$ and the part of the real axis contained in $H_\la$ belongs to~$W_{F_\la}$.
More generally, if $F$ is sufficiently close to $F_\la$, then the saddle fixed point $p_F$ of $F$ is contained in $W_F$ and the stable manifold $W^s(p_F) $ contains an arc close to the real axis connecting $p_F$ to the boundary of~$\DH$.
This arc is contained in $ W_F.$
Hence, every quasi-angular arc contained in $ \DH $ that turns once around the origin intersects $W_F$.

To every infinite backward orbit $ \underline{z} = \{z_m\}_{m \ge
0} \in \Gamma_F $ we associate a global unstable manifold $
W^u(\underline{z}) = \cup_{m \ge 0}
F^m(W^u_\alpha(\underline{z}_m))$, where $ \underline{z}_m $ is
the infinite backward orbit $ \{z_{m+j}\}_{j \ge 0}.$
 Below, in Lemma~\ref{wild}, we will prove  that there is $M \ge 1$ such that for every $\underline{z} \in \Gamma_F$ starting at a point in $\AFF$ and every $m = 1, \ldots, M$; we have that $F^m(W^u_\alpha(\underline{z})) \subset \DH$ and that $F^M(W^u_\alpha(\underline{z}))$ is a quasi-angular arc which turns once around the origin.
It follows that the local unstable manifold $W^u_\alpha(\underline{z})$ is contained in the global unstable manifold of an infinite backward orbit contained in $W_F$.
So this proves part~$1$ of Theorem~\ref{t-wild for endos}.

Recall that for $ F \in \sFla$ close to $ F_\la,$ the distance
from the fundamental domain $ \AFF $ to the origin has the order of
$ (1 - \la)^{-1} $, and the distance from the boundary of the
fundamental domain $ \DH $ to the origin  has  the order of $ (1 - \la)^{-\frac{1}{2}}.$ So, the distance from the fundamental domain $
\AFF $ of $ F$ to the boundary of the hyperbolic set $ \DH $ tends
to infinity as $ \lambda $ tends to $ 1.$

\begin{lemma}\label{wild}
For every $\lambda$ sufficiently close to~$1$ there exists an integer $ M \ge 0 $ such that for every $ F \in \sFla$ sufficiently
close to $ F_\la$, and every infinite backward orbit
$\underline{z} = \{ z_k \}_{k \ge 0}$ in $\Gamma_F$ starting at
$z_0 \in \AFF$, the following properties hold.

\begin{enumerate}

\item[1.]  For every $j = 0, \ldots, M$ the set
$F^j(W^u_\alpha(\underline{z})) $ is contained in the domain of
hyperbolicity $\DH$.

\item[2.] $F^M(W^u_\alpha(\underline{z})) $ is a quasi-angular arc that turns at least once around the origin.
\end{enumerate}
\end{lemma}

\proof
Let $L > 0$ be a constant such that every quasi-angular of length at least $L(1 - \la)^{-1}$ that is contained in $B_\la$, turns at least once around the origin.
On the other hand, choose $\eta \in (1, (5/2)^{\frac{1}{6}})$ and let $\la_1 \in (0, 1)$ be sufficiently close to~$1$ so that for every $\la \in (\la_1, 1)$ we have $\eta^3 < \la (5/2)^{\frac{1}{2}}$.

Notice that for every $w \in \C$ such that $|w| > R > 1$ we have $ |F_\la(w)| > \la R - 1 $.
Inductively, we have that $|F^m_\la(w)| > \la^m R - m$ whenever $ \la^{m - 1} R - (m - 1) > 1$.
Take $ R = \tfrac{1}{2}\rad (1 - \la)^{-1} < {\rm dist}(\{0\}, \AF_{F_\la})$ and consider an infinite backward orbit $ \underline{z} $ in $\Gamma_{F_\la}$ starting from a point $ z_0 $ in $ \AFla. $ Since for $ w \in
W^u_\alpha(\underline{z}) $ we have that $ |w| > R $, then
$F^m_\la(W^u_\alpha(\underline{z})) \subset \DH $ whenever $\la$ is sufficiently close to~$1$ and $\la^{m - 1}\tfrac{1}{2}\rad(1 - \la)^{-1} - (m - 1) > 5(1 - \la)^{-\frac{1}{2}}.$
Let $C > 0$ be such that for every $m \ge 0$ we have $C \eta^m (1 - \la)^{-\frac{1}{2}} \ge m - 1 + 5(1 - \la)^{-\frac{1}{2}}$.
Then, whenever
$$
m < \frac{\ln ((2C\la)^{-1} \rad) - \tfrac{1}{2}\ln (1 - \la)}{\ln (\eta/\la)},
$$
we have $F^m_\la(W^u_\alpha(\underline{z})) \subset \DH $.

On other hand, while the iterates of $ W^u_\alpha(\underline{z}) $
remain in $ \DH $ their length is increased by a factor $\la (5/2)^\frac{1}{2} > \eta^3 > 1$ (part~$1$ of Lemma~\ref{l-unstable cone field}).
So, the length of $ F^m_\la(W^u_\alpha(\underline{z})) $ is greater than $ 2 \eta^{3m} \alpha.$
Moreover, we have $ 2 \eta^{3m} \alpha > L (1 - \la)^{-1}$ if and only if
$$
m >
\frac{\ln (L / (2\alpha)) - \ln (1 - \la)}{3 \ln \eta}.
$$

Let $\la \in (\la_1, 1)$ be sufficiently close to~$1$ so that,
$$
\frac{\ln (L / (2\alpha)) - \ln (1 - \la)}{3 \ln \eta}
<
\frac{\ln ((2C\la)^{-1} \rad) - \tfrac{1}{2}\ln (1 - \la)}{\ln (\eta/\la)},
$$
and so that there is an integer $ M $ satisfying
$$
\frac{\ln (L / (2\alpha)) - \ln (1 - \la)}{3 \ln \eta}
< M <
\frac{\ln ((2C\la)^{-1} \rad) - \tfrac{1}{2}\ln (1 - \la)}{\ln (\eta/\la)}.
$$
Then the desired properties hold for $F = F_\la$ and for maps $F \in \sFla $ close enough to~$F_\la$.
\finproof

\subsection{A repelling annulus}
\begin{lemma}\label{repelling annulus lemma}
There is a closed annulus $\AFla^0 \subset \AFla$ that
is contained in the interior of its image under $F_\lambda$.
Moreover, $\AFla^0 $ contains at least one iterated preimage of
the fixed source $p^+_\la$ in its interior.
\end{lemma}

\proof Recall that the outer boundary $\gamma_{F_\lambda}^+$ of
$\AFla$ is a Jordan curve formed by a piece of the
unstable manifold of the saddle fixed point $p_\lambda$ that
intersects itself transversely at some point. Let
$\hgamma_\lambda^+$ be a Jordan curve in the interior of
$\AFla$ formed by a $C^1$ arc close to
$\gamma_{F_\lambda}^+$, that intersects itself transversely at
some point, and such that its image $F(\hgamma_\lambda^+)$ is
closer to the unstable manifold of $p_\lambda$, than
$\hgamma_\lambda^+$. Then the closed annulus $\AFla^0$ obtained
from $\AFla$ by replacing its outer boundary
$\gamma_{F_\lambda}^+$ by $\hgamma_\lambda^+$, is contained in
$\AFla$ and in the interior $F(\AFla^0)$. We can choose
$\hgamma_\lambda^+$ close enough to $\gamma_{F_\lambda}^+$ in such
a way that $\AFla^0$ contains an iterated preimage of
$p^+_\lambda$ in its interior. \finproof
\subsection{Curved arcs and tangencies}\label{curved arcs and tangencies}

Let $\AFla^0$ be the annulus given by Lemma~\ref{repelling annulus
lemma}. Denote by $\cC_\lambda$ the collection arcs $\gamma : [0,
1] \to \C $ of class $C^1$, such that the following properties
hold.
\begin{enumerate}
\item[1.] The image of $\gamma$ intersects $\AFla^0$ and $\length
(\gamma) \le d_\lambda = \dist (\AFla^0, \partial F_\lambda
(\AFla^0))/2$. \item[2.] For all $t \in [0, 1]$ the vector
$\gamma'(t)$ is non zero and it belongs to $\cK(\gamma(t))$.
\item[3.] There is $\rho \in \R$ (resp. $\rho' \in \R$) such that
vector $\gamma'(0)$ (resp. $\gamma'(1)$) is of the form
$$
\gamma'(0) = \gamma(0) \rho (i + 1/3) \  \ \left( \text{resp.
$\gamma'(1) = \gamma(1) \rho' (i - 1/3)$} \right).
$$
\end{enumerate}

Note that property~$1$ implies that $\gamma$ is contained in
$F_\lambda(\AFla^0)$ and that property~3 implies that $\gamma'(0)$
and $\gamma'(1)$ are in the boundary of $\cK(\gamma(0))$ and
$\cK(\gamma(1))$, respectively.

The following proposition is the key step to produce robust
tangencies.

\begin{proposition}\label{tangency to curved arc} If $\lambda$ is close to $1$ and $F$ in
$\sFla$ is close to $F_\lambda$, then every arc in
$\cC_\lambda$ is tangent to the local unstable manifold of an
infinite backward in~$\Gamma_F$.
\end{proposition}
\proof Consider $F \in \sFla$ close to $F_\lambda$ in such a
way that $\dist (\AFla^0, \partial F(\AFla^0)) > d_\lambda$.
Property~$1$ of the definition of $\cC_\lambda$ implies that the
image of every arc in $\cC_\lambda$ is contained in $F(\AFla^0)$.

\noindent {\bf 1.} We will show that for any arc $\gamma$ in
$\cC_\lambda$ there is a lift $\tgamma$ by $F$, such that a
sub-arc
 of $\tgamma$, re-parameterized to be defined in $[0, 1]$, belongs to $\cC_\lambda$.
In fact, as $\AFla^0 \subset F(\AFla^0)$, there is a lift
$\tgamma$ of $\gamma$ by $F$, such that for some $t \in [0, 1]$ we
have $\tgamma(t) \in \AFla^0$. Moreover, for each $t \in [0, 1]$
the vector $\tgamma'(t)$ belongs to $\tcK(\gamma(t))$, so
$\length(\tgamma) \le 2/3 \cdot \length(\gamma) < d_\lambda$.

On the other hand, since the set $F^{-1}(F(\AFla^0))$ is contained
in $\DH$, it follows that for every $z_0$ in this set we have
$D_{z_0}F(\cK(z_0)) \subset \cK(F(z_0))$. Hence, there are $\tau,
\tau' \in [0, 1]$ and $\rho, \rho' \in \R$, such that
$\tgamma'(\tau)$ and $\tgamma'(\tau')$ are of the form
$\gamma(\tau) \rho (i + 1/3)$ and $\gamma(\tau') \rho' (i - 1/3)$,
respectively. By taking $\tau'$ closer to $\tau$ if necessary, we
assume that for all $t$ between $\tau$ and $\tau'$, we have
$\tgamma'(t) \in \cK(\tgamma(t))$. So the arc $t \mapsto
\tgamma(\frac{t - \tau}{\tau' - \tau})$, defined on $[0, 1]$,
belongs to $\cC_\lambda$.

\noindent {\bf 2.} Let $\gamma_0$ be an arc in $\cC_\lambda$. For
$m \ge 1$ define inductively an arc $\gamma_m \in \cC_\lambda$ in
such a way that $\gamma_m$ is constructed from $\gamma_{m - 1}$,
as described in part~$1$. It follows that there is a sequence $\{
t_m \}_{m \ge 1} \subset [0, 1]$ such that for $m \ge 1$ we have
$F(\gamma_m(t_m)) = \gamma_{m - 1}(t_{m - 1})$, so that
$\underline{w} = \{ \gamma_m(t_m) \}_{m \le 0}$ is an infinite
backward orbit of points in $F_\lambda(\AFla^0) \subset \DH$. We
have then $\underline{w} \in \Gamma_F$. Moreover, the vector
$\gamma_m'(t_m)$ is non-zero and belongs $\cK(\gamma_m(t_m))$ and
for $m \ge 1$, the vector $DF(\gamma_m'(t_m))$ is parallel to
$\gamma_{m - 1}'(t_{m - 1})$. It follows that $\gamma_0'(t_0)$ is
parallel to the unstable direction associated to the backward
orbit $\underline{w}$ and that the arc $\gamma$ is tangent to the
corresponding local unstable manifold at $t = t_0$. \finproof

\subsection{The stable manifold of $p_F$ contains a curved arc}\label{the stable manifold is curved}
Note that every arc of length at most $d_\la$ which intersects $\AFla^0$ and
that turns around some point at least ~$2$ times contains an arc
in $\cC_\la$.

\noindent{\bf 1.} Define $R_0 = \{ t \in \C \mid t \in \R, t \le 0
\}$ and for $m \ge 1$ define inductively $R_m$ as the preimage of
$R_{m - 1}$ contained in $\H$, by the limit map $\Gcan$. The
corresponding backward orbit of ~0 converges to the fixed point
$p^+$ by Lemma~\ref{basin of p^+}.

\noindent{\bf 2.} As the derivative of $p^+$ has complex
multipliers, it follows that there is a positive integer $N $ and
a piece $\gamma$ of $R_N$ contained in $D$, such that $\gamma$ and
all its successive preimages under $\Gcan|_\H$ 'turn around $p^+$ at
least~$5$ times'.

\noindent{\bf 3.} As $F_\lambda$ converges in the $C^2$ topology
to $\Gcan$ on $D$, when $\lambda \rightarrow 1$, it follows that if
$\lambda$ is sufficiently close to~$1$, then there is a curve
$\gamma_\lambda$ contained in $D$ such that
$F_\lambda^N(\gamma_\la)$ is a piece of $R_0$, and such that
$\gamma_\lambda$ and all its successive preimages under
$F_\lambda|_\H$ 'turn around $p^+_{F_\lambda}$ at least~$4$
times'.

\noindent{\bf 4.} Let $z_0 \in A_\lambda^0$ be an iterated
preimage of $p^+_{F_\lambda}$ by $F_\lambda$. Then there is a
preimage $\tgamma_\lambda$ of $\gamma_\lambda$ by some iterate of
$F_\lambda$ such that $\tgamma_\lambda$ is contained in a
$d_\lambda / 2$ neighborhood of $z_0$, it has length at most
$d_\lambda /2$, and it 'turns around $z_0$ at least~$3$ times'.

\noindent{\bf 5.} So for $F \in \sFla$ sufficiently close to
$F_\lambda$ there is a piece of the stable manifold $\tgamma_F$
that is contained in a $d_\lambda$ neighborhood of $z_0$, it has
length at most $d_\lambda$ and it 'turns around $z_0$ at least~$2$
times'. It follows that a piece of $\tgamma_F$ is contained in
$\cC_\lambda$ and it is therefore tangent to an unstable manifold
of an infinite backward orbit in $\Gamma_F$.

\section{The vector fields $X_{\lambda,\mu}$}\label{vector field} This section contains the proof of
Theorem~\ref{embedding in R5} which consists of the construction
of a family of vector fields $\{ X_{\la, \mu} \}$.
In Subsection~\ref{remarks embedding} we state a refined
version of Theorem~\ref{embedding in R5},  make some remarks and
introduce  some notation. In Subsection~\ref{strategy embedding}
we describe the strategy of the construction of $X_{\la,\mu}$ and
give the  initial steps. The heart of the construction is
contained  in subsections~\ref{cherry flow} and~\ref{isotopy
section}.
We end the construction in Subsection~\ref{end construction}.
\subsection{Remarks and notation}\label{remarks embedding}
Recall that for an integer $k \ge 0$ we denote by $D^k$ the closed unit ball of $\R^k$.
Also $\overline{\D}$ is the closed unit disc in $\C$ and $\bfT^n = \R/\Z \times D^{n-1}$.

Fix an integer $n \ge 5$, $\la \in (0, 1)$ and $\eta > \sigma >
0$. Given $\mu \in (0, \sigma]$, let $F_{\la, \mu} : \C^* \to \C$
be defined by
$$
F_{\la, \mu}(z) = (1 - \la + \la |z|^{\mu/\sigma})(z/|z|)^2 + 1.
$$
Note that $F_{\la, \mu}$ depends smoothly on $\la$ and  $\mu$ and
that $F_{\la, \sigma} = F_\la$. For $ \lambda $ close
to~$1$ and $\mu \in (0, \sigma]$ close to $\sigma$,
the map $F_{\la, \mu}$ is area expanding.

Observe that $F_{\la, \mu}$ can be written as the composition
$F_{\la, \mu} = \Gcan \circ T_{\la, \mu}$ of the maps
$$\begin{array}{rccl}
T_{\la, \mu} : & \C^* &  \to & \C^*  \\
 & z &  \mapsto & (1 - \la + \la |z|^{\mu / \sigma}) ({z} / {|z|}),\\
\Gcan : &  \C^* & \to  & \C \\
& z & \mapsto & ({z^2}/{|z|}) + 1.
\end{array}$$

Recall that $B_\la = \{ z \in \C \mid |z| \le 2(1 - \la)^{-1} \}$, $B_\la^* = B_\la \setminus \{ 0 \}$ and put
$$
\tB_\la = \{ z \in \C \mid (1 - \la)/2 \le |z| \le 2(1 - \la)^{-1} - 1 - \la / 2 \} \subset B_\la.
$$
Then for all  $\mu \in (0,\sigma]$  we have
$T_{\la, \mu}(B_\la^*) \subset \tB_\la$, $\Gcan(\tB_\la) \subset B_\la$ and $F_\la(B_\la^*) \subset B_\la$.

\medskip

Similarly, for $\beta >0$ sufficiently small, the map
$$
\begin{array}{rccl}
  \hF_{\la,\mu} : & B^*_\la \times \overline{\D} \times D^{n-5} & \to         & B_\la \times \overline{\D} \times D^{n-5} \\
                         & (z,w,v)                                                         &\mapsto & \left( F_{\la, \mu}(z),  \dfrac{z}{ 2 |z|}
 + \beta |z|^{\eta/\sigma} \dfrac{|z|}{z} w ,  \beta |z|^{\eta/\sigma} v \right).
\end{array}
$$
is well defined and injective.
This map may be written as the
composition $\hF_{\la, \mu} = \hGcan \circ \hT_{\la, \mu}$ of the
maps
$$\begin{array}{rccl}
\hT_{\la, \mu} : & B^*_\la \times \overline{\D} \times D^{n - 5} & \to & \tB_{\la, \mu} \times \overline{\D} \times D^{n - 5} \\
& (z, w, v) & \mapsto & (T_{\la, \mu} (z), \beta_0 |z|^{\eta /
\sigma} w, \beta_0 |z|^{\eta / \sigma}v),
\end{array}$$
where $\beta_0 \in (0, 1)$ is sufficiently small, and
$$
\begin{array}{rccl}
\hGcan : & \tB_\la  \times \overline{\D} \times D^{n - 5} & \to & B_\la \times \overline{\D} \times D^{n - 5} \\
& ( z,   w,  v) & \mapsto & \left(\Gcan(z), \dfrac{z}{2|z|} + \beta_1 \dfrac{|z| }{z} w,
\beta_1  v\right),
\end{array}$$
where $\beta_1 \in (0, 1/2)$ is such that $\beta_0 \beta_1 =
\beta$.

\medskip
From the considerations above, Theorem~\ref{embedding in R5} is a
direct consequence of the following theorem.
\begin{theorem}\label{embedding in R5 II}
Fix an integer $n \ge 5$ and $\eta> \sigma >0$. Then for each
$\beta \in (0,1)$ sufficiently small and  for each $\la \in (0,1)$
sufficiently close to $1$,  there exists  $\mu_0  \in (0,\sigma)$
and a smooth one parameter family of smooth vector fields  $\{
X_{\la,\mu} \mid  \mu \in (\mu_0,\sigma] \}$ defined on an open
set $\trapp$ of  $\bfT^n$  such that for all $\mu \in (\mu_0,\sigma]$
the following hold:

\begin{enumerate}
\item[1.]
The boundary of the open set $\trapp  \subset \bfT^n$ is a manifold of dimension~$n - 1$ that is contained in the interior of $\bfT^n$.
For every $\mu \in [\mu_0, \sigma]$ the vector field $X_{\la, \mu}$ extends to a smooth vector field defined on a neighborhood of the closure of $\trapp$, in such a way that on the boundary of $\trapp$ this vector field points inward.
Moreover this extension has a unique singularity $o=o_{\la,\mu}$.
The singularity $o$ is contained in $\trapp$ and is hyperbolic with eigenvalues $-\mu$, $\sigma$ and $-\eta$ of multiplicities $1$, $2$ and $n-3$, respectively.

\item[2.]  There exist  codimension~$1$  submanifolds
$$
\Sigma^u \approx B_\la \times \overline{\D} \times D^{n - 5} \
\text{ and } \ \Sigma^s \approx \tB_\la  \times \overline{\D}
\times D^{n - 5},
$$
which are  transversal to the flow of $X_{\la,\mu}$ and so that
every forward orbit of the flow of $X_{\la,\mu}$ in $\trapp$
intersects $\Sigma^u$ or is contained in a local stable manifold
$W^s_{loc}(o)$ of $o$. The intersection  of this local stable
manifold $W^s_{loc}(o)$ with $\Sigma^u$ is $\{ 0 \} \times
\overline{\D} \times D^{n-5}$.

\item[3.] The Poincar\'e maps from $\Sigma^u$ to $\Sigma^s$ and
from $\Sigma^s$ to $\Sigma^u$ induced by $X_{\la, \mu}$ are given
by $\hT_{\la, \mu}$ and $\hGcan$, respectively.
\end{enumerate}
\end{theorem}

\subsection{Strategy of the construction}\label{strategy embedding}
In the rest of this section, we identify $\R^n$ with $\R \times \C
\times \C \times \R^{n-5}$ where we use coordinates   $(s,z,w,v)$.
Also, we view the solid torus $\bfT^n = \T \times D^{n-1}$ as a
subset of $\T \times \C \times \C \times \R^{n-5}$.

We will define a vector field $\tX_{\la,\mu}$ on a closed subset
$\tpretrapp$ of $[-2,3] \times 2B_\la \times 2\D \times D^{n-5}$
so that $\tX_{\la,\mu}$ is constant equal to $(1, {\bf 0})$ for
$s$ close to~$-2$ and~$3$. Hence, if we let $a=(1-\la)/16$, then
the image of $\tX_{\la,\mu}$  under the map
$$\begin{array}{rccc}
Q: & [-2,3] \times 2B_\la \times 2\D \times D^{n-5} & \rightarrow & \T \times D^{n-1} \\
   & (s,z,w,v)  & \mapsto & ( (s+2)/5 \mod 1,  a z,  a  w, v/4)
\end{array} $$
will be a smooth vector field $DQ(\tX_{\la,\mu})$ defined on the
subset $\pretrapp = Q(\tpretrapp)$ of $\T \times \frac{1}{4} \overline{\D}
\times \frac{1}{4} \overline{\D} \times \frac{1}{4} D^{n-5}$.  Our
desired vector field  $X_{\la,\mu}$ will be obtained extending
 $DQ(\tX_{\la,\mu})$ to a small neighborhood $\trapp$ of $\pretrapp$.

The construction of $\tX_{\la,\mu}$ will be such that the codimension~$1$ submanifolds $\tSigma^u_- = \{ -2 \} \times B_\la \times \overline{\D} \times
D^{n-5}$, $\tSigma^s = \{ 2 \} \times \tB_\la \times \overline{\D}
\times D^{n-5}$  and  $\tSigma^u_+ = \{ 3 \} \times B_\la \times
\overline{\D} \times D^{n-5}$ are transversal to $\tX_{\la, \mu}$ (See Figure~\ref{xlamu figure}).
The Poincar\'e map from
$\tSigma^u_-$ into $\tSigma^s$ will be given by the map
$(-2,z,w,v) \mapsto (2, \hT_{\la,\mu} (z,w,v))$. This Poincar\'e
map corresponds to a transition through a singularity with
eigenvalues $-\mu$ in the $s$-direction, $\sigma$ in
$z$-directions, and $-\eta$ in the $w$ and $v$ directions. The
Poincar\'e map from $\tSigma^s$ into $\tSigma^u_+ $ will be given
by $(2, z,w,v) \mapsto (3, \hGcan(z,w,v))$ and corresponds to an
isotopy between the identity and $\hGcan$.

The heart of the proof of Theorem~\ref{embedding in R5 II} is to
construct the vector field $\tX_{\la,\mu}$. In
Subsection~\ref{cherry flow} we construct the part corresponding
to the passage through a singularity and in
Subsection~\ref{isotopy section} we construct the part
corresponding to the isotopy. Then, in Subsection~\ref{end
construction} we finish the construction of $X_{\la,\mu}$.


\subsection{Through the singularity}\label{cherry flow}
Recall that we identify  $\R^n$ with $ \R \times \C \times \C
\times \R^{n-5}$ and use coordinates   $(s,z,w,v)$. In $\R^n$ we
consider the linear diagonal vector field $L_\mu$ which has
eigenvalue $-\mu$ in the $s$-direction, $\sigma$ in the two
$z$-directions, and $-\eta$ in all the $w$ and $v$ directions.
Although one should think of $L_\mu$ as a vector field defined on
a copy of $\R^n$ distinct from the one where we construct
$\tX_{\la,\mu}$, at one step of the construction will be convenient to
think of $L_\mu$ as defined in a neighborhood of $s=-1$ in the
copy of $\R^n$ which corresponds to the vector field $\tX_{\la,
\mu}$ .

\subsubsection*{Definition of  $\tX_{\la,\mu}$ for $s \in [-2,-1]$}
Let $\tX_{\la,\mu}$ be constant equal to $(1,{\bf 0})$ on a
neighborhood of $s=-2$ and let $\tX_{\la,\mu}$ coincide with $L_\mu$ on a neighborhood of $s=-1$.
Then, letting $b = 4(1 - \lambda)^{-1}$, extend $\tX_{\la,\mu}$ to $s \in (-2, -1)$ so that the Poincar\'e map from $s=-2$ to
$s=-1$ is
$$
(-2,z,w,v) \mapsto (-1, z/b , w, v).
$$
Thus the Poincar\'e map  shrinks $B_\la$ to the
closed disc  $\frac{1}{2}\overline{\D}$ of radius $1/2$ in the
$z$-coordinate.

Let $\tpretrapp_{[-2,-1]}$ be the set of points $(s,z,w,v)$ with $s
\in[-2,-1]$ which belong to an orbit of the flow of
$\tX_{\la,\mu}$ that starts at  $\tSigma^u_- = \{-2 \} \times
B_\la  \times \overline{\D} \times D^{n-5}$.

\subsubsection*{Definition of  $\tX_{\la,\mu}$ for $s \in [-1,1]$}
\begin{figure}
\begin{center}
\psfig{file=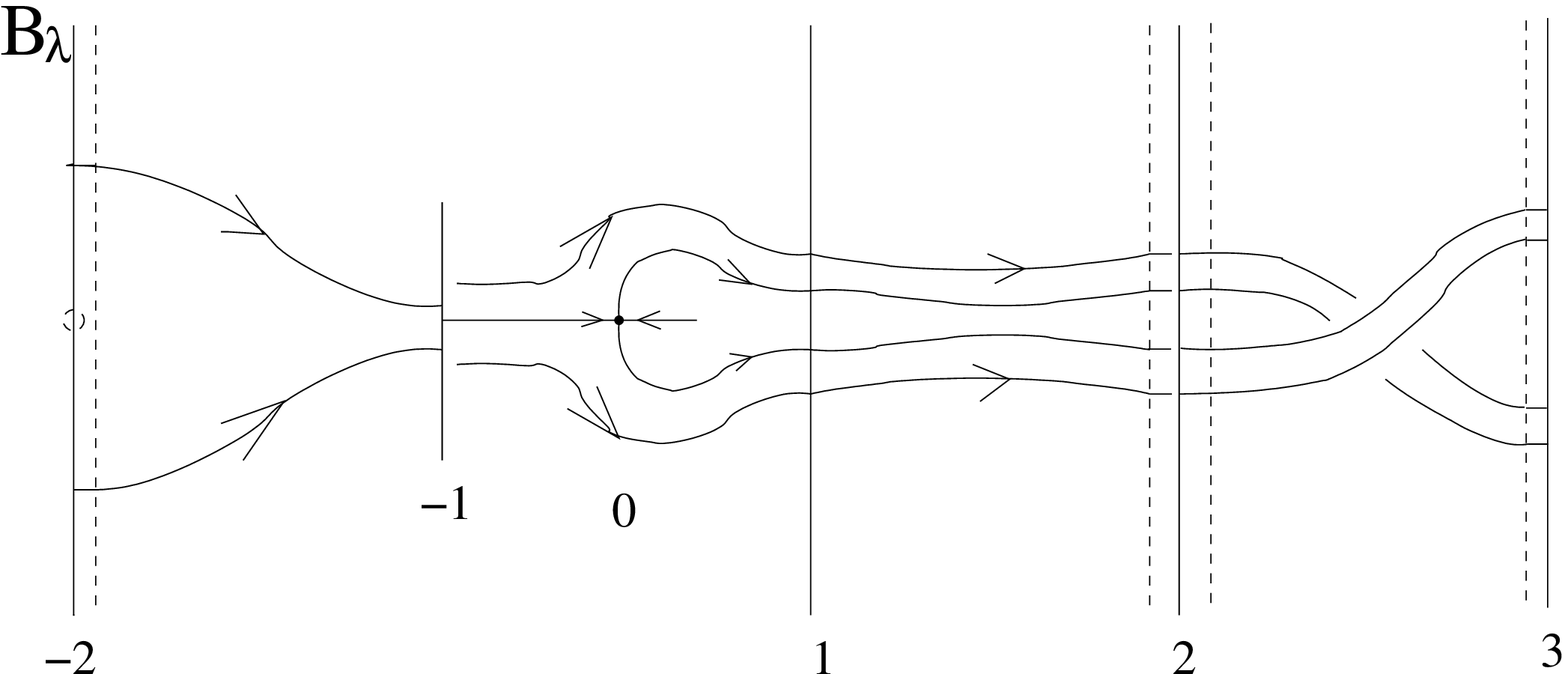,width=8cm}
\end{center}
\caption{Illustrates the construction of  $\tX_{\la,\mu}$} \label{xlamu figure}
\end{figure}

Now observe that, under the linear flow $L_\mu$ we have the
following  transition map through the singularity:
$$\begin{array}{ccl}
\{-1\} \times \frac{1}{2}\overline{\D}^*  \times \overline{\D} \times D^{n-5} & \rightarrow & [-2^{-\mu/\sigma}, 0) \times \{ |z|=1\} \times \overline{\D} \times D^{n-5} \\
(-1,z,w,v) & \mapsto & (-|z|^{\mu/\sigma}, z/|z|,
|z|^{\eta/\sigma} w, |z|^{\eta/\sigma}v).
\end{array}$$
For $\epsilon >0$ small, let $\linearpiece$ be the closed region in $\R^n$
bounded by $s=-1$, $s=\epsilon$, $|z| = 1$ and the orbits of
$\{-1, \epsilon\} \times \partial \left( \frac{1}{2} \overline{\D}
\right) \times \overline{\D} \times D^{n-5}$  under the linear
flow $L_\mu$ (see Figure~\ref{xlamu figure}).

Choose $\beta_0 >0$ small and consider a smooth family of maps
$$\psi_\mu : \linearpiece   \rightarrow [-1,1] \times 2B_\la \times 2\D \times D^{n-5}$$
such that the following hold:
\begin{enumerate}
  \item[(a)] $\psi_\mu$ extends to a smooth family of  diffeomorphisms between a neighborhood of $\linearpiece$
and a neighborhood $\psi_\mu(\linearpiece)$.
  \item[(b)]  For all  $(s,z,w,v)$ in a neighborhood of $\{-1\} \times \frac{1}{2}\overline{\D}  \times \overline{\D}  \times D^{n-5}$,
 $$\psi_\mu (s,z,w,v) = (s,z,w,v).$$
\item[(c)] For all $(s,z,w,v) \in [-2^{-\mu/\sigma},
2^{-\mu/\sigma} \epsilon] \times \{ |z|=1\} \times \overline{\D}
\times D^{n-5}$,
$$\psi_\mu (s,z,w,v) = (1, (1-\la+ \la b^{\mu/\sigma}(-s)) z, \beta_0 b^{\eta/\sigma} w,\beta_0 b^{\eta/\sigma} v)$$
\end{enumerate}

For example we may consider a non-increasing $C^\infty$ function
$h_0:[0,1] \rightarrow [0,1]$ such that $h_0(r) = 0$ for $r \in [0, 1/2]$, $h_0$ is strictly increasing on $[1/2,1]$ and $h_0(1) = 1$.
Then we may take $\psi_\mu (s,z,w,v) = (s',z',w',v')$ where,
\begin{eqnarray*}
  s' & =& h_0(|z|)(1 - s) + s \\
  z' &=& \left( h_0(|z|)(- \la +\la b^{\mu/\sigma}(-s)) + 1 \right)z \\
  w' & =& \left( h_0(|z|)(\beta_0 b^{\eta/\sigma} -1) +1 \right) w \\
  v' &= &  \left( h_0(|z|)(\beta_0 b^{\eta/\sigma} -1) +1 \right) v.
\end{eqnarray*}
A standard exercise shows that such $\psi_\mu$ is a diffeomorphism onto its image.

Now define   $\tX_{\la,\mu}$ on $\tpretrapp_{[-1,1]} = \psi_\mu(\linearpiece)$ as $D
\psi_\mu (L_\mu)$. Note that the Poincar\'e map between
$\tSigma^u_-$  and $s=1$ is given by
$$(-2,z,w,v) \mapsto (1,\hT_{\la,\mu} (z,w,v)).$$
Denote by $\tsing$ the unique singularity $\psi_\mu ({\bf 0})$ of
$\tX_{\la,\mu}$ in $\tpretrapp_{[-1,1]}$.

\subsubsection*{Definition of  $\tX_{\la,\mu}$ for  $s \in [1,2]$}
We let $\tX_{\la,\mu}$ be constant equal to $(1,{\bf 0})$ on a
neighborhood of $s=2$ and extend $\tX_{\la,\mu}$ to $s \in (1,2)$ so
that the Poincar\'e map from $s=1$ to $s=2$ is the identity in the
$(z,w,v)$ coordinates.

Let $\tpretrapp_{[1,2]}$ be the set of points $(s,z,w,v)$ with $s \in
[1,2]$ which belong to an orbit starting at $\psi_\mu (\linearpiece) \cap \{
s= 1 \}$.

\subsection{The isotopy}\label{isotopy section}
Our aim now is to define $\tX_{\la, \mu}$ for $s \in [2,3]$ so
that $\tX_{\la,\mu}$ is the horizontal vector field $(1, {\bf 0})$
in a neighborhood of $s=2$ and $s=3$ and so that the Poincar\'e map
from $s=2$ to $s=3$ is given by $(2, z,w,v) \mapsto (3,
\hGcan(z,w,v))$. The idea is to construct an isotopy between the
identity and $\hGcan$.

We need the following Lemmas.

\begin{lemma}
Consider the map $ \gamma:[0, 1] \times \T \rightarrow \C^2 $
defined by
$$
\gamma_s(\theta) = ((1 - s) \exp({2\pi i \theta}) + s \exp({4 \pi i
\theta}), s \exp({2 \pi i \theta})).
$$
Then the following properties hold.
\begin{enumerate}
\item[1.]
For every $ s \in [0, 1]$ the map $\gamma_s$ is injective, and for every $\theta$ we have $\gamma'_s(\theta) \neq 0$.
In particular $ \gamma $ defines an isotopy between $ \gamma_0 $ and $ \gamma_1.$

\item[2.] There is a smooth map $u:[0, 1] \times \T \rightarrow
\C^2 $ such that $u_0(\theta) \equiv (0, 1)$, $u_1(\theta) = ( 0,
\exp({-2 \pi i \theta}))$ for every $\theta$, and for every $ s $ and
$ \theta$, the vectors $ u_s(\theta) $ and $ \gamma'_s(\theta) $
are linearly independent over $\C$.
\end{enumerate}

\end{lemma}
\proof The first assertion is easily proved. To prove the second
let $\chi:[0, 1] \rightarrow [0, 1]$ be a smooth function such
that $\chi \equiv 0$ on $ [0, 1/3]$ and $\chi \equiv 1$ on $[2/3,
1]$. For $0 < \varepsilon < 1/2$ small define,

$$
u_s(\theta) = \left\{
\begin{array}{lcl}
  (- 4 \chi(\varepsilon^{-1}s), 1 - \chi(\varepsilon^{-1}s)) & {\rm if} & s \in [0,\varepsilon] \\
  (-4, 0) & \text{if} & x \in [\varepsilon, 1 - \varepsilon] \\
  (- 4 (1 - \chi(1 - \varepsilon^{-1}(1 - s))),  \chi(1 - \varepsilon^{-1}(1 - s)) \exp({- 2 \pi i \theta})) & \text{if} & s \in [1 - \varepsilon, 1] \\
\end{array}
\right.
$$

The lemma is now proved by straightforward calculations. \finproof

\begin{lemma}
For $s \in [0,1]$ and $\beta_1 \in (0,1)$, let $\hG_s: \tB_\la \times \overline{\D} \times
D^{n-5} \to 2 B_\la \times 2 \D \times D^{n-5}$ be the map defined by
$$
\hG_s (z,  w,  v) = (\gamma_s(\theta) + \varepsilon (- i t
\gamma'_s(\theta) + w u_s(\theta)), \beta_1 v) 
$$
where $t = |z|$ and $  \theta = \tfrac{1}{ 2\pi} \arg(z) \in \T$.
If $ \beta_1 \in (0, 1) $ is sufficiently small, then for every $
s \in [0, 1] $ the map $\hG_s: \tB_\la \times \overline{\D} \times
D^{n-5} \to 2 B_\la \times 2 \D \times D^{n-5}$ 
is injective.
\end{lemma}
\proof Observe that  $\|\partial_\theta \hG_s -
(\gamma'_s(\theta), 0)\| \leq \varepsilon(M \|\gamma''_s(\theta)
\| + \| u'_s(\theta)\|)\;\;\; {\rm and} \;\;\;
\partial_t \hG_s = (- \varepsilon i \gamma'_s(\theta), 0).$
Moreover, $ \hG_s $ is holomorphic in $ w $ and $ \partial_w \hG_s
= (\varepsilon u_s(\theta), 0).$ So, if $\varepsilon $ is
sufficiently small, the vectors $\partial_\theta \hG_s $ and $
\partial_t \hG_s$ are linearly independent over $ \R$, and the plane
generated by real combinations of $\partial_\theta \hG_s $ and $
\partial_t \hG_s$ is close to the complex plane $\{ (\la
\gamma'_s(\theta), 0)  \mid \la \in \C \}.$

By construction $\gamma'_s(\theta) $ and $ u_s(\theta) $ are
linearly independent over $ \C.$ It follows that, if $\varepsilon
$ is sufficiently small, $\hG_s $ is a local diffeomorphism, and
that there is $ \delta > 0 $ such that for every $ \theta_0 $ the
map $ \hG_s $ restricted to $\{(t \exp({2 \pi i \theta}), w, v) \mid
|\theta - \theta_0| < \delta \} $ is injective. As $ \gamma_s $ is
injective for all $ s \in [0, 1]$, it follows that, if $
\varepsilon $ is sufficiently small, $ \hG_s $ is injective for
all $ s.$ \finproof

Clearly, $\hG_0(z, w, v) = ((z / |z|) + 2 \pi \varepsilon z,
\varepsilon w, \beta_1 v) $ is isotopic to the identity, and
$$\hG_1(z, w, v) = \left({z^2 \over |z|^2} + 4 \pi \varepsilon {z^2 \over |z|},
 {z \over |z|} + 2 \pi \varepsilon z + \varepsilon w {|z| \over z}, \beta_1 v \right) $$ is
isotopic to $ \hGcan.$

\

Let $\widehat H: \tB_\la \times \overline{\D} \times D^{n-5} \to 2
B_\la \times 2 \D \times D^{n-5} $ be an isotopy between the
identity and $\hGcan$. Let $h_1 : [2, 3] \to [0, 1]$ be a
$C^{\infty}$ function that is constant equal to~$0$ (resp.~$1$) on
a neighborhood of~$2$ (resp.~$3$). The vector field $
\tX_{\la,\mu} $ will be defined on the set
$$
\tpretrapp_{[2,3]} = \left\{ \left( s, \widehat H_{h_1(s)}({\bf
x}) \right) \mid  s \in [2, 3] \ {\rm and } \,\, {\bf x} \in
\tB_\la \times \overline{\D} \times D^{n-5} \right\},
$$
by
$$
 \tX_{\la,\mu} (s,   z,   w,   v)
= \left(1, \frac{\partial}{\partial s} \left( \widehat H_{h_1(s)}
\left( \widehat H_{h_1(s)}^{-1}(z,w,v) \right) \right) \right).
$$

By construction $ \tX_{\la,\mu}$ induces the Poincar{\'e} map $(2,
{\bf x}) \mapsto (3, \hGcan({\bf x}))$. Since the function $h_1$ is
constant on a neighborhood of $s = 2$ and of $s = 3$, it follows that $
\tX_{\la,\mu}$ is constant equal to $(1, {\bf 0})$ on a
neighborhood of $s=2$ and of $s=3$.

\subsection{End of the construction}\label{end construction}
We let  the domain of definition of the vector field
$\tX_{\la,\mu}$ be
$$\tpretrapp = \tpretrapp_{[-2,-1]} \cup \tpretrapp_{[-1,1]} \cup \tpretrapp_{[1,2]} \cup \tpretrapp_{[2,3]}$$
and observe that the unique singularity $\tsing$ of $\tX_{\la,\mu}$
is $\psi_\mu({\bf 0}) \in \tpretrapp_{[-1,1]}$. The stable manifold of this
singularity $W^s(\tsing)$ consists of $\psi_\mu (\linearpiece \cap \{ z = 0
\}) \subset \tpretrapp_{[-1,1]}$ together with the points in $\tpretrapp_{[-2,-1]}$
that are in the backward orbit of points in  $\psi_\mu (\linearpiece \cap \{
z = 0 \})$. In particular, $W^s(\tsing) \cap \tSigma^u_- = \{ -1 \}
\times \{0 \} \times \overline{\D} \times D^{n-5}$. Moreover, for
every neighborhood $N$ of $W^s(\tsing)$, there exists $t_0 >0$ such
that the orbit of every point $(s, {\bf x}) \in V \setminus N$
hits $\tSigma^u_+$ before time $t_0$. Note that $W^s(\tsing) \cap
\tSigma^u_+ =\emptyset$ and $W^s(\tsing) \cap \tSigma^s =\emptyset$

Now we pass to the solid torus $\bfT^n$ by the map $Q: (s,z,w,v)
\mapsto ((s+2)/5, a z,  a  w,$ $ v/4)$, and obtain the vector
field $X_{\la,\mu} = DQ (\tX_{\la,\mu})$ defined on  $\pretrapp =
Q(\tpretrapp)$ whose  flow is transversal to the cross sections $\Sigma^u
= Q(\tSigma^u_\pm)$ and $\Sigma^s = Q(\tSigma^s)$. Moreover,
$X_{\la,\mu}$ has a unique singularity $o= Q(\tsing)$ in $\pretrapp$ with
local stable manifold $W^s_{loc}(o) = Q(W^s (\tsing))$.

Take $\mu = \sigma$ and observe that for every $t_0 > 0$ we have $X^{t_0}_{\la,\sigma} (\pretrapp)
\subset \operatorname{interior}(\pretrapp)$. Therefore, we may extend the
definition of $X_{\la, \sigma}$ to a neighborhood $\pretrapp'$ of $\pretrapp$ and find a neighborhood $\trapp$ of $\pretrapp$ that is bounded by a smooth manifold of dimension $n - 1$ contained in $\pretrapp'$, in such a way that $X_{\la, \sigma}$ points inward on the boundary of $\trapp$.
Since $\trapp$ contains $\pretrapp$ for $\mu$ close to $\sigma$ we may extend the
definition of $X_{\la,\mu}$ from $\pretrapp$ to a neighborhood of the closure of $\trapp$, so that we have a
smooth family of vector fields defined in $\trapp$ for $\mu$ close to
$\sigma$. It follows that there exists $\mu_0$ so that for all
$\mu \in (\mu_0, \sigma]$ the vector field $X_{\la, \mu}$ points inward on the boundary of $\trapp$.
The rest of the properties required in Theorem~\ref{embedding in R5 II} for  $X_{\la, \mu}$ easily
follow.

\section{First return map and leaf space transformation}\label{s-dimension reduction}
In this section we show that for vector fields $X$ close to $X_{\la, \mu}$, the first return map $\hF_X$ to $\Sigma^u$ of the flow of $X$ admits a strong stable foliation of codimension~$2$ (Subsection~\ref{ss-first return map}) and we study the corresponding leaf space transformation (Subsection~\ref{ss-leaf space transformation}).
In particular, we show that  when $X$ is of class $C^2$, the corresponding leaf space transformation belongs to $\sFla$ (Lemma~\ref{l-leaf space transformation}).
\subsection{First return map}\label{ss-first return map}

Fix an integer $n \ge 5$. For a given $\la \in (0, 1)$ close
to~$1$, let $X_{\la} = X_{\la, \sigma}$ be the vector field
defined on ${\bf T}^n$ and let $\Sigma^u \approx B_\la \times
\overline{\D} \times D^{n - 5}$ be the transversal section to
$X_\la$, given by Theorem~\ref{embedding in R5 II}.

Note that the intersection between a local stable manifold of the
singularity $\sing = \sing_{\la, \sigma}$ and $\Sigma^u$ is equal to $\{ 0 \}
\times \overline{\D} \times D^{n - 5}$.
For $X$ close to $X_\la$ in the $C^1$ topology we denote by $\sing_X$ the singularity of $X$ that is the continuation of the hyperbolic singularity $\sing$ of $X_\la$.
After a smooth coordinate change, we assume that for every vector field $X$
close to $X_\la$ in the $C^1$ topology, the intersection of a local stable manifold of the singularity $o_X$ of $X$ with $\Sigma^u$ is equal to $\{ 0 \}
\times \overline{\D} \times D^{n - 5}$.
Then we set $\Sigma^{u*} = \Sigma^u \setminus (\{ 0 \} \times \ov{\D} \times D^{n - 5})$ and observe that there is a well defined Poincar\'e map $\hF_X :
\Sigma^{u*} \to \Sigma^u$ of the flow of~$X$.

The following lemma is a consequence of~\cite{HPS}.
Recall that $\P1 : \Sigma^u \approx B_\la \times \ov{\D} \times D^{n - 5} \to B_\la$ is the projection to the first coordinate.
\begin{lemma}[Strong stable foliation]\label{l-ss foliation}
Let $\la \in (0, 1)$ be close to~$1$, let $\mu \in (0, \sigma)$
be close to $\sigma$ and let $X_{\la, \mu}$ be the vector field given
by Theorem~\ref{embedding in R5 II}. Then there is a neighborhood
$\cO$ of $X_{\la, \mu}$ in the $C^1$ topology, such that for every
$X$ in $\cO$ we have the following properties.
\begin{enumerate}
\item[1.] There is a strong stable foliation $\folss_X$ of $\hF_X$ in
$\Sigma^u$, having $\{ 0 \} \times \ov{\D} \times D^{n - 5}$ as a leaf.
The leaves of $\folss_X$ are
submanifolds of codimension 2 in  $\Sigma^u$ that are of class
$C^1$. Moreover the foliation $\folss_X$ is $C^0$ close to the
foliation formed by the fibers of the projection~$\P1$.
 \item[2.]
Every leaf of the strong stable foliation intersects $B_\la \times
\{ 0 \} \times \{ {\bf 0} \}$ in at most one point. For every
point $\bfx \in \Sigma^u$ in a leaf intersecting this set, we
denote by $\Pi_X(\bfx)$ the point in $B_\la$ such that
$(\Pi_X(\bfx), 0, {\bf 0})$ is in the same leaf as $\bfx$. Then
the map $\Pi_X$ is continuous and $C^0$ close to the projection
$\P1$. In the particular case when $X$ is of class $C^2$, the map
$\Pi_X$ is of class $C^1$ and it depends on a $C^1$ way on $X$.
\end{enumerate}
\end{lemma}
\proof As $\mu \in (0, \sigma)$, it follows that $\hF_{\la, \mu}$
contracts the fibers of the projection $\P1$ in a stronger way
than any contraction of $F_{\la, \mu}$.

Moreover, the stronger contraction factor of $F_{\la, \mu}$ in
$B_\la^*$,
$$
\inf\{ \| D(F_{\la, \mu})_z (v) \| \mid z \in B_\la^*, |v| = 1 \},
$$
can be made arbitrarily close to $ \la$ by choosing $\mu
\in (0, \sigma)$ sufficiently close to ~$\sigma$. On the other
hand, the fiber contraction can be made arbitrarily strong by
choosing $ \beta > 0 $ small enough. Then the results follow from
~\cite{HPS} using a graph transforming method. For a direct
exposition of these methods which applies to our case see
~\cite{BLMP}. \finproof

\subsection{Leaf space transformation}\label{ss-leaf space transformation}
Let $\la \in (0, 1)$, $\mu \in (0, \sigma)$ and $\cO$ be as in Lemma~\ref{l-ss foliation}. As the closure of $\hF_{\la, \mu}(\Sigma^{u*})$ is in the
interior of $\Sigma^u$, reducing $\cO$ if necessary we assume that
every leaf through a point in $\hF_X(\Sigma^{u*})$ intersects $B_\la
\times \{ 0 \} \times \{ {\bf 0} \}$ in a unique point. Then the
{\it leaf space  transformation} $F_X : B_\la^* \to B_\la$ is defined
by
$$
F_X(z) = \Pi_X(\hF_X(z, 0, {\bf 0})).
$$
The map $F_X$ is continuous, but in general not differentiable.
\begin{lemma}\label{l-local homeo}
Reducing $\cO$ if necessary we have that for every $X \in \cO$ the map $F_X$ is a local homeomorphism.
\end{lemma}

\proof
Note that two points $z, z' \in \C^*$ have the same image under $F_{\la, \mu}$ if and only if $z + z' = 0$.
On the other hand, the images under $\hF_{\la, \mu}$ of the leaves $\{z\} \times \D \times D^{n-5}$ and $\{-z\} \times \D \times D^{n-5}$ of $\folss_{X_{\la, \mu}}$ lie in the same leaf of $\folss_X$ and their distance is at least $2(1 - \beta).$

By continuity, for every vector field $X$ that is sufficiently close to~$X_{\la, \sigma}$ in the $C^1$ topology, the distance between~$2$ points in the same leaf of $\folss_X$ whose preimages by $\hF_X$ are in distinct leaves, is at least $ 1 - \beta.$
This implies that arbitrarily close
leaves of $ \folss_X $ can not have images contained in the  same leaf.
It follows that $F_X$ is locally injective and therefore a local homeomorphism.
\finproof

For $C^2 $ vector fields we can say even more.

\begin{lemma}\label{l-leaf space transformation}
Denote by $\cO'$ the subset of $\cO$ of vector fields of class
$C^2$. Then, reducing $\cO$ if necessary, the following properties hold.
\begin{enumerate}
\item[1.]
For each vector field $X$ in $\cO'$, there are $\tau_X \in \sTla$ and $g_X \in \sGla$ such that $F_X = g_X \circ \tau_X$.
In particular $F_X \in \sFla$.
\item[2.]
The maps $X \mapsto \tau_X$, $X \mapsto g_X$ and $X \mapsto F_X$ from $\cO'$ to $\sTla$, $\sGla$ and $\sFla$, respectively, are all continuous. Note that here $\cO' \subset \cO$ is provided with the $C^1$ topology.
\item[3.]
For each $X \in \cO'$ there is a neighborhood $\sU$ of $g_X$ in $\sGla$, such that for every $g \in \sU$ there is $Y \in \cO$ satisfying $F_Y = g \circ \tau_X$.
\item[4.]
For each $X \in \cO'$ the map $F_X$ is area expanding.
\end{enumerate}
\end{lemma}

\proof

\partn{1}
For $ X \in \cO' $ the first return map $ \widehat F_X $ is of
class $ C^2$ and $\Pi_X$ is of class $C^1$.
So $F_X$ is of class $C^1$.

For each vector field $X$ that is close to $X_\la$ in the $C^1$ topology, let $\hT_X$ (resp. $\hG_X $) be the Poincar\'e map from $ \Sigma^u $  to $\Sigma^s$ (resp. from $ \Sigma^s $ to $  \Sigma^u$).
Clearly $\hF_X = \hG_X \circ \hT_X$.
The pull-back  of the strong stable foliation $\folss_X$ from $
\Sigma^u $  to $ \Sigma^s $ by $\hG_X $ defines a foliation $\widetilde{\folss_X}$ in $ \Sigma^s \approx \tB_\la  \times \overline{\D} \times D^{n - 5}$.
The map $\hT_X $ carries leaves of $\folss_X$ into leaves of $\widetilde{\folss_X}$.
When $X = X_\la$, the foliation $\widetilde{\folss_{X_\la}}$ is the one formed by the fibers of the projection $\tP1 : \tB_\la \times \overline{\D} \times D^{n - 5} \to \tB_\la$ to the first coordinate.
Thus the foliation $\widetilde{\folss_X}$ is close to $\widetilde{\folss_{X_\la}}$.
Let $\tPi_X$ be the projection along leaves of $\folss_X$, in such a way that $\bfx$ and $(\tPi_X(\bfx), 0, \bf0)$ are in the same leaf of $\widetilde{\folss_X}$.
The map $\tPi_X$ is defined on a large subset of $\tB_\la \times \overline{\D} \times D^{n - 5}$ and takes images in $\tB_\la$.
Note that $\tPi_{X_{\la}} = \tP1$ and that when $X$ of class $C^2$ varies continuously in the $C^1$ topology, $\tPi_X$ varies in a $C^1$ way.

Let $\iota : B_\la \to B_\la \times \ov{\D} \times D^{n - 5}$ (resp. $\tiota : \tB_\la \to \tB_\la \times \ov{\D} \times D^{n - 5}$) be defined by $z \mapsto (z, 0, {\bf 0})$ and put
$$
T_X = \tPi_X \circ \hT_X \circ \iota : B_\la^* \to \tB_\la,
$$
$$
G_X = \Pi_X \circ \hG_X \circ \tiota : \tB_\la \to B_\la.
$$
Note that $F_X = G_X \circ T_X$, $\hT_{X_\la} = \hT_{\la, \sigma}$, $\hG_{X_\la} = \hGcan$, $T_{X_\la} = T_{\la, \sigma}$ and $G_{X_\la} = \Gcan$.
So if we denote by $\psi : \C^* \to \T \times \R$ the map $\psi(z) = (\tfrac{1}{2\pi}\arg z, |z|)$, with inverse $\psi^{-1}(\theta, t) = t \exp(2\pi i \theta)$, then $\tau_\lambda = \psi \circ T_{X_\la}$ and $\gcan = G_{X_\la} \circ \psi^{-1}$.

For a vector field $X$ of class $C^2$ that is close to $X_\la$ in the $C^1$ topology, the image of $B_\la^*$ by $T_X$ is an annulus that is bounded by the Jordan curve $T_X(\{ z \in \C \mid |z|  = 2(1 - \la)^{-1}\})$ and the Jordan curve equal to the image by $\tPi_X$ of the intersection of the local unstable manifold of $\sing_X$ with $\Sigma^s$.
These Jordan curves are of class $C^1$ and they vary continuously in the~$C^1$ topology, when $X$, of class $C^2$, varies continuously in the $C^1$ topology.
So for such $X$ we can find a homeomorphism,
$$
\psi_X : \overline{T_X(B_\la^*)} \to \T \times [1 - \la, 2(1 - \la) - 1 - \la],
$$
that extends to a diffeomorphism onto its image, defined on a neighborhood of $\overline{T_X(B_\la^*)}$, in such a way that $\psi_{X_\la}$ coincides with $\psi$ and that $\psi_X$ varies continuously in the (strong) $C^1$ topology when $X$, of class $C^2$, varies continuously in the $C^1$ topology.
Then we have $\tau_X = \psi_X \circ T_X \in \sTla$, $g_X = G_X \circ \psi_X^{-1} \in \sGla$ and $F_X = G_X \circ T_X = g_X \circ \tau_X \in \sFla$.

\partn{2}
Clearly when $X$ of class $C^2$ varies continuously in the $C^1$ topology, $G_X$, $\psi_X$ and $g_X = G_X \circ \psi_X^{-1}$ vary continuously in the (strong) $C^1$ topology.
It remains to show that $T_X$ varies continuously in the (weak) $C^1$ topology.
For that, just observe that for each neighborhood $\hU$ of $\{0\} \times \bar \D \times D^{n-5}$ in $\Sigma^u$ (the intersection of the stable manifold of $o_{\la, \mu}$ with $ \Sigma^u$) the map $X \mapsto \hT_X
\mid_{(\Sigma^u \setminus \hU)}$ is continuous in the $C^1$ topology.

\partn{3}
Keep the notation of part~$1$.
Note first that for every map $\hG$ of class $C^1$ that is close to $\hG_X$ there is a vector field $Y \in \cO$ that coincides with $X$ between $\Sigma^u$ and $\Sigma^s$ (so that $\hT_Y = \hT_X$) and such that $\hG_Y = \hG$.
On the other hand, note that for every $g$ close to $g_X$ we can find $\hG$ close to $\hG_X$ mapping leaves of $\widetilde{\folss_X}$ into leaves of $\folss_X$ and such that $\Pi_X \circ \hG \circ \tiota$ coincides with $g \circ \psi_X$ on $\tau_X(B_\la^*)$.

For a given $g \in \sGla$ near $g_X$ let $\hG$ and $Y$ be as above, so that $\hT_Y = \hT_X$, $\hG_Y = \hG$ and $\hG_Y$ maps leaves of $\widetilde{\folss_X}$ into leaves of $\folss_X$.
It follows that the strong stable foliation $\folss_Y$ of $\hF_Y$ is equal to $\folss_X$ and that $\Pi_Y = \Pi_X$.
So $G_Y = \Pi_X \circ \hG \circ \tiota$ coincides with $g \circ \psi_X$ on $\tau_X(B_\la^*)$ and we have,
$$
F_Y = G_Y \circ T_Y = g \circ (\psi_X \circ T_X) = g \circ \tau_X.
$$

\partn{4}
To prove that $ F_X $ is area expanding, we first notice that the map
$F_{\la, \mu} $ is area expanding.
For each vector field $ X \in \cO'$ the map $F_X $ is $C^1$-close to $F_{\la, \mu}
$ outside a small neighborhood of the origin. So $F_X$ is area
expanding on this set.
That $F_X$ is area expanding close to~$0$ is clear from the eigenvalues of the linear part of $ X $ at the singularity $\sing_X$, as these are close to those of the linear part of $X_{\la, \mu} $ at $\sing_{\la, \mu}$. \finproof


\section{Robust transitivity}\label{s-robust transitivity}
This section is dedicated to the proof of Theorem~\ref{t-dynamics of the flow}.
\subsection{Proof of part~$1$ of Theorem~\ref{t-dynamics of the flow}}\label{ss-part 1 C2}
Let $X_{\la, \mu}, \trapp, \sing, \ldots$ be as in the statement of Theorem~\ref{embedding in R5}.
For a vector field $X$ close to $X_{\la, \mu}$, let $\Lambda_X$ be the maximal invariant set in $\trapp$ of the flow of $X$ and let $o_X$ be the singularity of $X$ that is the continuation of the hyperbolic singularity~$\sing$ of~$X_{\la, \mu}$.
Clearly we have $\sing_X \in \Lambda_X$.
By construction of $X_{\la, \mu}$, it follows that for every vector field~$X$ that is sufficiently close to $X_{\la, \mu}$, the point~$\sing_X$ is the unique singularity of $X$ in $\trapp$.
So, to prove part~1 of Theorem~\ref{t-dynamics of the flow} it remains to show that for $\la \in (0, 1)$ sufficiently close to~$1$ and $\mu \in (0, \sigma)$ sufficiently close to $\sigma$, there is a neighborhood $\cO$ of $X_{\la, \mu}$ in the $C^1$ topology such that for every $X \in \cO$ the restriction of the flow of $X$ to $\Lambda_X$ is topologically mixing.

Recall that $\hF_X : \Sigma^{u*} \to \Sigma^u$ is the first return map to $\Sigma^u$ of the flow of $X$.
Although $\hF_X$ is not defined on $\{ 0 \} \times \ov{\D} \times D^{n - 5}$ we let it act on subsets of $\Sigma^u$ by $\hF_X(\hU) := \hF_X(\hU \setminus \{ 0 \} \times \ov{\D} \times D^{n - 5})$.
We denote by
$$
\hOmega_X = \cap_{m \ge 1} \hF_X^m(\Sigma^{u}),
$$
the maximal invariant set of $\hF_X$ in $\Sigma^u$.
The maximal invariant set $\Lambda_X$ of the flow of~$X$ is just the closure of the suspension of $\widehat{\Omega}_X$ by the flow of~$X$.
So to prove part~$1$ of Theorem~\ref{t-dynamics of the flow} is enough to prove that the restriction of $\hF_X$ to~$\hOmega_X$ is topologically mixing.

In view of Lemma~\ref{l-leaf space transformation}, part~$1$ of Theorem~\ref{t-dynamics of the flow} is a direct consequence of the following proposition and of Theorem~\ref{t-topological dynamics}.
Note that for every $X \in \cO$ the maximal invariant set $\Omega_X = \Omega_{F_X}$ of $F_X$ in $B_\la$ is equal to $\Pi_X(\hOmega_X)$.
\begin{proposition}\label{p-part 1 C2}
Let $\la \in (0, 1)$, $\mu \in (0, \sigma)$ and $\cO$ be as in Lemma~\ref{l-leaf space transformation}.
Then for every vector field $X \in \cO$ such that the restiction of $F_X$ to $\Omega_{X}$ is topologically mixing, the map $\hF_X$ restricted to $\hOmega_X$ is topologically mixing.
\end{proposition}
For the proof of Proposition~\ref{p-part 1 C2} we need the following lemma.
\begin{lemma}
Given an open set $\hU$ in $ \Sigma^u  $ intersecting $\hOmega_X$,
there exists an open set $ U_0 $ in $B_\la$ intersecting $
\Omega_X $ and an integer $j \ge 0$, such that
$\hF_X^j(\Pi_X^{-1}(U_0)) \subset \hU.$
\end{lemma}
\proof Let $ \hU $ be an open set as in the hypothesis and take $
\bfx \in \hU  \cap \hOmega_X. $ Then for every integer $j \ge 0$
there exists a point $\bfx_j \in \hOmega_X $ such that $
\widehat{F}_X^j(\bfx_j) = \bfx. $ As $\hF_X$ contracts the leaves
of the foliation $\folss_X$ uniformly, for~$j$ sufficiently
large the leaf of $\folss_X$ through the point $\bfx_j $ is
mapped well inside $ \hU $ by $ \hF_X^j.$ The same happens for a
neighborhood of leaves of $\folss_X$. \finproof

\proofof{Proposition~\ref{p-part 1 C2}}
Let $\hU, \hV$ be two open subsets of $\Sigma^u$ that intersect $ \hOmega_X$.
Take an open set $U_0 $ in $ B_\la $ and a positive integer $j$ such that
$\hF_X^j (\Pi_X^{-1}(U_0)) \subset \hU$, given by the previous
lemma. Since the restriction of $F_X$ to $\Omega_X$ is topologically mixing and $U_0$ intersects $\Omega_X$, there exists a positive
integer $k_0$ such that for every $k \ge k_0$ we have
$F^k_X(\Pi_X(\hV)) \cap U_0 \neq \emptyset$. Hence for every $k
\ge k_0$ we have $\hF_X^k (\hV) \cap \Pi_X^{-1}(U_0) \neq
\emptyset$. But this implies that for every $\ell \ge j + k_0$ we
have $\hF^\ell_X (\hV) \cap \hU \neq \emptyset $. So, the restriction of $\hF_X $ to $\hOmega_X$ is topologically mixing.
\finproof

\subsection{Proof of part~$2$ of Theorem~\ref{t-dynamics of the flow}}\label{ss-part 2}
Let $\la \in (0, 1)$, $\mu \in (0, \sigma)$ and $\cO$ be as in Lemma~\ref{l-leaf space transformation}.
Let $\hp_\la$ be the unique fixed point of $\hF_{\la, \sigma}$ in $\P1^{-1}(p_\la)$.
Reducing $\cO$ if necessary, we assume that for each $X \in \cO$ there is a fixed point $\hp_X$ of $\hF_X$ that is the continuation of $\hp_\la$ and such that $\Pi_X(\hp_X) \in \DH$.
Thus, for every $X \in \cO$ of class $C^2$ we have $\Pi_X(\hp_X) = p_{F_X}$.

For $X \in \cO$ let
$$
\hW_X = \{ \bfx \in \hOmega_X \mid \hF_X^m(\bfx) \in \Pi_X^{-1}(\DH) \text{ for all $m \ge 0$} \},
$$
be the maximal invariant set of $\hF_X$ in $\Pi_X^{-1}(\DH)$.
Clearly $\hW_X$ is a compact set invariant by $\hF_X$.

\partn{1}
Identify the tangent space at each point of
$\Sigma^u \approx B_\la \times \overline{\D} \times D^{n - 5}$
with $\C \times \C \times \R^{n - 5}$ and denote by $\| \cdot \|$
the Euclidean norm in $\R^{n - 5}$. Put $\varepsilon_0 =
\tfrac{1}{2}(1 - \la)^{1/2}$, as in Subsection~\ref{ss-stable cone field}, and let $\hcC$ and $\hcK$ be the cone fields defined on
$B_\la^* \times \overline{\D} \times D^{n - 5}$ by
\begin{multline*}
\hcC(z_0, v_0, w_0) = \left\{ (z_0\rho(1 + i\varepsilon), \mathbf{w},\mathbf{v})
\mid \rho \ge 0, \ |\varepsilon| \le \varepsilon_0 \right.
\\
\left. \text{ and } \ \rho \ge (1 - \la)\max \{ |\mathbf{w}|, \| \mathbf{v} \| \}
\right\},
\end{multline*}
\begin{multline*}
\hcK(z_0, v_0, w_0) = \left\{ (z_0\rho(i + \varepsilon), \mathbf{w}, \mathbf{v})
\mid \rho \ge 0, \ |\varepsilon| \le 1/ 3 \right.
\\
\left. \text{ and } \ \rho \ge (1 - \la)\max \{ |\mathbf{w}|, \| \mathbf{v} \| \}
\right\}.
\end{multline*}
With a straightforward computation is easy to check that, if
$\beta > 0$ in the definition of $F_{\la, \mu}$ is sufficiently
small, then for every $\la \in (0, 1)$ sufficiently close to~$1$
and $\mu \in (0, \sigma)$ sufficiently close to $\sigma$, the cone
field $\hcC$ (resp. $\hcK$) is an stable (resp. unstable) and
invariant cone field for $\hF_{\la, \mu} : \P1^{-1}(\DH) \to
\hF_{\la, \mu}(\P1^{-1}(\DH))$.
Then, reducing $\cO$ if necessary, for every $X \in \cO$ the
cone field $\hcC$ (resp. $\hcK$) is an stable (resp. unstable) and
invariant cone field for $\hF_X : \Pi_X^{-1}(\DH) \to
\hF_X(\Pi_X^{-1}(\DH))$.
It follows that $\hW_X$ is a uniformly hyperbolic set for $\hF_X$.

\partn{2}
For $X \in \cO$ put,
$$
\hGamma_X = \{ \{ \bfx_j \}_{j \ge 0} \mid \Pi_X(\bfx_j) \in \DH,
\ \hF_X(\bfx_{j + 1}) = \bfx_j \}.
$$
Note that for every infinite backward orbit $\{ \bfx_j \}_{j \ge
0}$ of $\hF_X$ in $\hGamma_X$, the sequence $\{ \Pi_X(\bfx_j)
\}_{j \ge 0}$ is an infinite backward orbit of $F_X$ in $\Gamma_X
:= \Gamma_{F_X}$. Conversely, for every infinite backward orbit
$\{ z_j \}_{j \ge 0}$ of $F_X$, the sequence $\{ \bfx_j \}_{j \ge
0}$ defined by the property $$ \{ \bfx_{j_0} \} = \cap_{j \ge 0}
\hF_X^j(\Pi_X^{-1}(z_{j_0 + j})),$$ is an infinite backward orbit of
$\hF_X$ in $\hGamma_X$.
From the theory of uniformly hyperbolic sets we know that for every infinite backward orbit $\widehat{\underline{z}} \in \hGamma_X$ there is a one dimensional local unstable manifold
$W_\alpha^u(\widehat{\underline{z}})$ of $\hF_X$ through $\bfx_0$.

Reducing $\cO$ if necessary, part~1 of Theorem~\ref{t-wild for endos} implies that, if $X \in \cO$ is of class $C^2$, then the local unstable manifold of each infinite backward orbit in $\hGamma_X$ is contained in the unstable manifold of some point of $\hW_X$.

\partn{3} Let $\la \in (0, 1)$ be sufficiently close
to~$1$ and let $\sU$ be a sufficiently small neighborhood of
$F_\la$ in $\sFla$, so that for every $F \in \sU$ there is an
arc $\tgamma_F$ of the stable manifold of the fixed point $p_F$ of
$F$ that is tangent to the local unstable manifold of an infinite backward orbit $\underline{z}_F$ of $F$ in $\Gamma_F$ (part~2 of Theorem~\ref{t-wild for endos}).

Reducing $\cO$ if necessary we assume that for every $X \in \cO$ of class $C^2$ we have $F_X \in \sU$ (Lemma~\ref{l-leaf space transformation}).
As all the elements of $\underline{z}_{F_X}$ belong to $\AF_{F_X}$, which
is well inside $\DH$, it follows that there is an infinite
backward orbit $\widehat{\underline{z}}_X$ of $\hF_X$ in $\hGamma_X$
that projects to $\underline{z}_{F_X}$ by $\Pi_X$. Therefore the local
unstable manifold $W_\alpha^u(\widehat{\underline{z}}_X)$ is tangent
to the submanifold $\Pi_X^{-1}(\tgamma_{F_X})$ of the stable
manifold of~$\widehat{p}_X$.

\noindent{\bf 3.} Let $\la \in (0, 1)$ be sufficiently close
to~$1$ and let $\sU$ be a sufficiently small neighborhood of
$F_\la$ in $\sFla$, so that for every $F \in \sU$ there is an
arc $\tgamma_F$ of the stable manifold of the fixed point $p_F$ of
$F$ that is tangent to the local unstable manifold of an infinite backward orbit $\underline{z}_F$ of $F$ in $\Gamma_F$ (part~2 of Theorem~\ref{t-wild for endos}).

Reducing $\cO$ if necessary we assume that for every $X \in \cO$ of class $C^2$ we have $F_X \in \sU$ (Lemma~\ref{l-leaf space transformation}).
As all the elements of $\underline{z}_{F_X}$ belong to $\AF_{F_X}$, which
is well inside $\DH$, it follows that there is an infinite
backward orbit $\widehat{\underline{z}}_X$ of $\hF_X$ in $\hGamma_X$
that projects to $\underline{z}_{F_X}$ by $\Pi_X$. Therefore the local
unstable manifold $W_\alpha^u(\widehat{\underline{z}}_X)$ is tangent
to the submanifold $\Pi_X^{-1}(\tgamma_{F_X})$ of the stable
manifold of~$\widehat{p}_X$.

\partn{4}
Let $X$ be an arbitrary vector field in $\cO$
and let $\{ X_j \}_{j \ge 0}$ be a sequence of vector fields of
class $C^2$ in $\cO$, that converge to $X$ in the $C^1$ topology.
Let $\underline{\bfx}^j = \widehat{\underline{z}}_{X_j} \in \hGamma_{X_j}$ be given by part~$3$, so that the local unstable manifold of $\underline{\bfx}^j$ is tangent to $\Pi_X^{-1}(\tgamma_{F_X})$.
Setting $\underline{\bfx}^j = \{ \bfx_k^j \}_{k \ge 0}$ and taking a subsequence if necessary, we
assume that for every $k \ge 0$ the $\bfx_k^j$ converge to some
point $\bfx_k$ as $j \to \infty$. It follows that $\Pi_X(\bfx_k)
\in \DH$ and that $\underline{\bfx} = \{ \bfx_k \}_{k \ge
0}$ is an infinite backward orbit of $\hF_X$ in $\hGamma_X$.
Moreover is easy to check that the local unstable manifolds
$W_\alpha^u(\underline{\bfx}^j)$ converge to
$W_\alpha^u(\underline{\bfx})$ in the $C^1$ topology and
that the manifolds $\Pi^{-1}_{X_j}(\tgamma_{F_{X_j}})$ converge in the
$C^1$ topology to a submanifold of the stable manifold of
$\widehat{p}_{X}$, as $j \to \infty$. So
$W^u(\underline{\bfx})$ is tangent to the stable manifold
of $\widehat{p}_X$.
As $\widehat{p}_X$ is contained in $\hW_X$ and by part~$2$ the set $W^u(\underline{\bfx})$ is contained in the unstable foliation of $\hW_X$, it follows that $\hW_X$ is a wild hyperbolic set for $\hF_X$.

\subsection{Proof of part~$3$ of Theorem~\ref{t-dynamics of the flow}}\label{ss-part 3}
Let $\sU$ be the neighborhood of $F_\la$ in $\sFla$ given by Proposition~\ref{p-periodic points}.
Let $\la \in (0, 1)$ sufficiently close to~$1$, $\mu \in (0, \sigma)$ sufficiently close to $\sigma$ and $\cO$ a sufficiently small neighborhood of $X_{\la, \mu}$ such that for every $X \in \cO$ of class $C^2$ we have $F_X \in \sU$ (Lemma~\ref{l-leaf space transformation}).

\partn{1}
We will show that there is a dense subset $\cD$ of $\cO$ of vector fields $X$ such that the set of periodic sources and the set of periodic saddles of $F_X$ are both dense in~$\Omega_X$.
To prove this assertion, let $Y \in \cO$ be a given vector field of class $C^2$ and let $\tau_Y \in \sT$ and $g_Y \in \sGla$ be given by Lemma~\ref{l-leaf space transformation}.
Part~$4$ of the same lemma implies that $g_Y \in \sU_{\tau_Y}$ and Proposition~\ref{p-periodic points} implies that there is $g$ close to $g_Y$ in $\sGla$ such that the set of periodic sources and the set of periodic saddles of $F = g \circ \tau_Y$ are both dense in $\Omega_F$.
By part~$3$ of Lemma~\ref{l-leaf space transformation} there is a vector field $X \in \cO$ near $Y$ such that $F_X = g \circ \tau_Y$.
This shows the assertion.

\partn{2}
We will now prove that for $X \in \cD$ the set
of periodic points of Morse index~$1$ (resp.~$2$) of $\hF_X$ is dense
in $\hOmega_X$.

Given a point $\bfx$ in $\hOmega_X$ and an integer $j \ge 0$ let $\bfx_j$ be
the unique point in $\Sigma^u$ such that $\hF_X^j(\bfx_j) = \bfx$.
Let $\{ p_j \}_{j \ge 0}$ be a sequence of saddle (resp.
repelling) periodic points of $F_X$ close enough to
$\Pi_X(\bfx_j)$, so that $F_X^j(p_j)$ converges to $\Pi_X(\bfx)$
as $j \to \infty$. Since $\hF_X$ contracts uniformly the fibers of
the map $\Pi_X$, it follows that the sequence of periodic point of
Morse index~$1$ (resp. $2$) $\{ \hF_X^j(\widehat{p}_j) \}_{j \ge 0}$
of $\hF_X$ converges to $\bfx$.

\partn{3}
Let $\hT$ be a trapping region for
$\hF_{X_{\la, \mu}}$ containing $\hOmega_{F_{\la, \mu}}$. Reducing
$\cO$ if necessary we assume that for every $X$ in $\cO$, the set
$\hT$ contains $\hOmega_X$ and is a trapping region for $\hF_X$.
Let $\hcB$ be a countable base of the topology of $\hT$.

For $ \hU \in \hcB $  and $ X \in \cO \cap \cD $ we have either $
\widehat{\Omega}_X \cap \hU = \emptyset $ or  $ \widehat{\Omega}_X
\cap \hU \ne \emptyset. $ In the former case there exists a $C^1$
open neighborhood $ \cO_X^{\hU} $ of $ X $  such that for every $Y \in \cO_X^{\hU}$ we have $\widehat{\Omega}_Y \cap \hU = \emptyset$.
In the later case there exists a $C^1$ open neighborhood $\cO_X^{\hU}$ of $ X $
 such that for every $ Y \in \cO_X^{\hU}$ the set $ \hU $ contains a periodic saddle (resp. source) of $\hF_Y$.

Then, $ \cD_{\hU} = \cup_{X \in \cO \cap \cD} \cO_X^{\hU} $ is an $C^1$
open and dense subset of $\cO$ and $ \cR = \cap_{\hU \in \hcB} \cD_{\hU}
$ is a $C^1$ residual subset of $\cO$ such that for every $ Y \in
\cR$, the set of periodic points of Morse index~$1$ (resp.~$2$) of
$\hF_Y$ is dense in $ \hOmega_X.$
\finproof

\bibliographystyle{plain}

\end{document}